\documentclass[11pt,a4paper]{article}
\usepackage{amsmath, amssymb}
\usepackage{latexsym, color}
\usepackage{epsfig}

\usepackage{pgf,tikz}
\usetikzlibrary{arrows}

\numberwithin{equation}{section}

\begin{document}

\newtheorem{thm}{Theorem}[section]
\newtheorem{cor}[thm]{Corollary}
\newtheorem{lem}[thm]{Lemma}
\newtheorem{prop}[thm]{Proposition}
\newtheorem{definition}[thm]{Definition}
\newtheorem{rem}[thm]{Remark}
\newtheorem{Ex}[thm]{Example}
\def\nm{\noalign{\medskip}}

\bibliographystyle{plain}


\newcommand{\qed}{\hfill \ensuremath{\square}}
\newcommand{\ds}{\displaystyle}
\newcommand{\pf}{\noindent {\sl Proof}. \ }
\newcommand{\p}{\partial}
\newcommand{\pd}[2]{\frac {\p #1}{\p #2}}
\newcommand{\norm}[1]{\| #1 \|}
\newcommand{\dbar}{\overline \p}
\newcommand{\eqnref}[1]{(\ref {#1})}
\newcommand{\na}{\nabla}
\newcommand{\one}[1]{#1^{(1)}}
\newcommand{\two}[1]{#1^{(2)}}

\newcommand{\Abb}{\mathbb{A}}
\newcommand{\Cbb}{\mathbb{C}}
\newcommand{\Ibb}{\mathbb{I}}
\newcommand{\Nbb}{\mathbb{N}}
\newcommand{\Kbb}{\mathbb{K}}
\newcommand{\Rbb}{\mathbb{R}}
\newcommand{\Sbb}{\mathbb{S}}

\renewcommand{\div}{\mbox{div}~}

\newcommand{\la}{\langle}
\newcommand{\ra}{\rangle}

\newcommand{\Hcal}{\mathcal{H}}
\newcommand{\Lcal}{\mathcal{L}}
\newcommand{\Kcal}{\mathcal{K}}
\newcommand{\Dcal}{\mathcal{D}}
\newcommand{\Pcal}{\mathcal{P}}
\newcommand{\Qcal}{\mathcal{Q}}
\newcommand{\Scal}{\mathcal{S}}

\def\Ba{{\bf a}}
\def\Bb{{\bf b}}
\def\Bc{{\bf c}}
\def\Bd{{\bf d}}
\def\Be{{\bf e}}
\def\Bf{{\bf f}}
\def\Bg{{\bf g}}
\def\Bh{{\bf h}}
\def\Bi{{\bf i}}
\def\Bj{{\bf j}}
\def\Bk{{\bf k}}
\def\Bl{{\bf l}}
\def\Bm{{\bf m}}
\def\Bn{{\bf n}}
\def\Bo{{\bf o}}
\def\Bp{{\bf p}}
\def\Bq{{\bf q}}
\def\Br{{\bf r}}
\def\Bs{{\bf s}}
\def\Bt{{\bf t}}
\def\Bu{{\bf u}}
\def\Bv{{\bf v}}
\def\Bw{{\bf w}}
\def\Bx{{\bf x}}
\def\By{{\bf y}}
\def\Bz{{\bf z}}
\def\BA{{\bf A}}
\def\BB{{\bf B}}
\def\BC{{\bf C}}
\def\BD{{\bf D}}
\def\BE{{\bf E}}
\def\BF{{\bf F}}
\def\BG{{\bf G}}
\def\BH{{\bf H}}
\def\BI{{\bf I}}
\def\BJ{{\bf J}}
\def\BK{{\bf K}}
\def\BL{{\bf L}}
\def\BM{{\bf M}}
\def\BN{{\bf N}}
\def\BO{{\bf O}}
\def\BP{{\bf P}}
\def\BQ{{\bf Q}}
\def\BR{{\bf R}}
\def\BS{{\bf S}}
\def\BT{{\bf T}}
\def\BU{{\bf U}}
\def\BV{{\bf V}}
\def\BW{{\bf W}}
\def\BX{{\bf X}}
\def\BY{{\bf Y}}
\def\BZ{{\bf Z}}


\newcommand{\Ga}{\alpha}
\newcommand{\Gb}{\beta}
\newcommand{\Gd}{\delta}
\newcommand{\Ge}{\epsilon}
\newcommand{\Gve}{\varepsilon}
\newcommand{\Gf}{\phi}
\newcommand{\Gvf}{\varphi}
\newcommand{\Gg}{\gamma}
\newcommand{\Gc}{\chi}
\newcommand{\Gi}{\iota}
\newcommand{\Gk}{\kappa}
\newcommand{\Gvk}{\varkappa}
\newcommand{\Gl}{\lambda}
\newcommand{\Gn}{\eta}
\newcommand{\Gm}{\mu}
\newcommand{\Gv}{\nu}
\newcommand{\Gp}{\pi}
\newcommand{\Gt}{\theta}
\newcommand{\Gvt}{\vartheta}
\newcommand{\Gr}{\rho}
\newcommand{\Gvr}{\varrho}
\newcommand{\Gs}{\sigma}
\newcommand{\Gvs}{\varsigma}
\newcommand{\Gj}{\tau}
\newcommand{\Gu}{\upsilon}
\newcommand{\Go}{\omega}
\newcommand{\Gx}{\xi}
\newcommand{\Gy}{\psi}
\newcommand{\Gz}{\zeta}
\newcommand{\GD}{\Delta}
\newcommand{\GF}{\Phi}
\newcommand{\GG}{\Gamma}
\newcommand{\GL}{\Lambda}
\newcommand{\GP}{\Pi}
\newcommand{\GT}{\Theta}
\newcommand{\GS}{\Sigma}
\newcommand{\GU}{\Upsilon}
\newcommand{\GO}{\Omega}
\newcommand{\GX}{\Xi}
\newcommand{\GY}{\Psi}

\newcommand{\beq}{\begin{equation}}
\newcommand{\eeq}{\end{equation}}

\title{Two types of electric field enhancements by\\ infinitely many circular conductors\\ arranged closely in two parallel line}

\author{KiHyun Yun \\{\small Department of Mathematics,}\\ {\small Hankuk University of Foreign Studies}\\{\small kihyun.yun@gmail.com}}

\maketitle

\begin{abstract}
 In this paper, we consider very high concentration of electric field in between infinitely many circular perfect conductors arranged closely in two rows. In stiff fiber-reinforced composite, shear stress concentrations occur in between neighboring fibers, and the electric field means shear stress in this paper. Due to material failure of composites, there have been intensive studies so far to estimate the field in between only a finite number of inclusions. Indeed, fiber-reinforced composites contain a large number of stiff fibers, and the concentration can be strongly enhanced by some combinations of inclusions. Thus, we establish some asymptotes and  optimal blow-up rates for the field in narrow regions in between infinitely many conductors in two rows to describe the effects combined horizontally and vertically by a large number of inclusions. Especially, the one of  blow-up rates is substantially different from the existing result in the case of finite inclusions.

\end{abstract}

\noindent {\footnotesize {\bf Mathematics subject classification
(MSC2000)}: 35J25, 73C40}

\noindent {\footnotesize {\bf Keywords}: conductivity equation, anti-plane elasticity, stress, blow-up, extreme conductivity}

\section{Introduction}
\par In stiff fiber-reinforced composites, high shear stress concentrations occur in between closely spaced neighboring fibers \cite{BC}. In the anti-plane shear model, the out-of-plane displacement $u$ satisfies a conductivity equation whose inclusions in the plane are the cross-sections of fibers, and the gradient $\nabla u$ implies the shear stress tensor. The problem to estimate $\nabla u$  in between inclusions was raised by Babu\u{s}ka  in the study of material failure of composites \cite{bab}. Many studies on the gradient estimate have been successfully carried out due to such practical significance \cite{LV,LN, LLBY,bv}. The genetic blow-up rate of $|\nabla u|$ is $\frac 1 {\sqrt {\tilde \epsilon}}$ for small $\tilde \epsilon >0$ when $\tilde \epsilon$ is the distance between two neighboring inclusions \cite{AKLim,AKLLL,Y,Y2}, and moreover, asymptotic behavior of $\nabla u$ was also established \cite{ACKLY,KLeeY,KLY}. The two dimensional problem has been generalized in various ways including high dimensions \cite{BLY,BLY2,KLY,KLY13,lekner10,lekner11,lekner,LLBY,LY,LYu,Y3}. Especially, it has been shown in \cite{LY2} that the concentration of $\nabla u$ can be strongly enhanced by a small inclusion between inclusions. This means that some combinations of inclusions can have strong influence on the concentration. So far, such studies have considered the cases when only a finite number of inclusions exist. 
 
 \par This paper is mainly concerned with the concentration of $\nabla u$ enhanced by a combination of infinitely many inclusions, because  composites contain a large number of stiff fibers. Thus, we consider infinitely many circular inclusions arranged closely in two rows to describe the horizontal and vertical effects of infinitely many inclusions. According to our results, one of effects is very strong enough to provide the blow-up rate substantially different from the existing rate in the case of finite number of inclusions.

\par We set up infinitely many circular perfect conductors arranged closely in two rows.  For any integer number $n$, we choose a pair of unit open disks $D_{Rn}$ and $D_{Ln}$ spaced $\epsilon$ apart in the horizontal direction, and moreover the distances between $D_{Rn}$ and $D_{Rn+1}$, and between $D_{Ln}$ and $D_{Ln+1}$ both are $\delta$ in the vertical direction. The open disks $D_{Rn}$ and $D_{Ln}$ are defined as $$ D_{R n} = \left\{ (x,y) | \left(x -\left (1 + \frac {\epsilon} 2\right) \right)^2 +  \left(y - n \left(2+\delta\right)\right)^2 < 1 \right\}$$
and
$$ D_{L n} = \left\{ (x,y) | \left(x - \left( -1 - \frac {\epsilon} 2\right) \right)^2 +  \left(y - n \left(2+\delta\right)\right)^2 < 1 \right\}.$$ Then, the domin $\mathbb{R}^2 \setminus \overline {\cup_{n=-\infty}^{\infty}\left(D_{R n} \cup D_{L n}\right)}$ has a periodic structure with period $2+\delta$ in the $y$ direction. In this paper, we suppose that $\epsilon$ and  $\delta $ are sufficiently small and positive.

\par Dealing with the governing equation, let the symbol $H$ denote a harmonic function defined in $\mathbb{R}^2$ whose gradient is a periodic function with period $2+\delta$ in the $y$ direction satisfying
\beq \nabla H(x,y) =  \nabla H(x,y + 2 + \delta) \mbox { for any } (x,y)\in \mathbb{R}^2. \label{H_condi}\eeq For example, a linear function $$H(x,y)= ax + by$$ can be a harmonic function with a periodic gradient described above. For a  such harmonic function $H$, we estimate the gradient $\nabla u$ of a solution  $u$  to the equation 
\beq  \Delta u = 0 \mbox{ in }  \mathbb{R}^2 \setminus \overline {\cup_{n=-\infty}^{\infty}\left(D_{L n} \cup D_{R n}\right)} \label{main}\eeq with the conditions \beq
 \begin{cases}
 \ds u = {c}_L   &\mbox{ on }\p D_{L 0},  \\ \nm   u = c_R    &\mbox{ on } \p D_{R 0}  ,\\ \nm \ds \int_{\p D_{L 0}} \p_{\nu} u ds = \ds \int_{\p D_{R 0}} \p_{\nu} u ds =0,&\\
  \nm \ds  \int_{\Omega \setminus \overline { D_{L0} \cup D_{R0}}} |\nabla (u-H)|^2 dxdy < \infty,  & \\  \nm \ds  \nabla u(x,y) =  \nabla u(x,y + 2 + \delta) & \mbox { for } (x,y)\in \mathbb{R}^2 \setminus \overline {\cup_{n=-\infty}^{\infty}\left(D_{L n} \cup D_{R n}\right)} .  \end{cases}\label{main_condi}
 \eeq Here, the constants $c_L$ and $c_R$ depend on $H$, $\epsilon$ and $\delta$, and  the normal unit vector $\nu$ points toward the inside of $D_{L0}$ or  $D_{R0}$. The domain $\Omega$ denotes a horizontal area as \beq \Omega = \mathbb{R} \times \left(-1 - \frac {\delta}  2,~ 1 + \frac {\delta}  2\right), \label{def_Omega}\eeq and the domain $\Omega_m$ is defined as \beq \Omega_m =\left \{ (x,y) ~|~(x,y)\in \Omega\mbox { and } |x|< m  \right\} \label{Omega_m}\eeq for $m \geq 3$ containing $D_{L0}$ and  $D_{R0}$.

   By definition, the gradient $\nabla u $ is a periodic function with period $2+\delta$ in the $y$ direction, and the solution $u$ has  a constant Dirichlet boundary data on each of boundaries $\p D_{R n} $ and $\p D_{L n}$ for  $n = 0,\pm 1, \pm 2, \cdots$. The existence of the solution $u$ for  a  harmonic function $H$ can be shown by considering $u(x,y)+u(x,-y)$ and $u(x,y) -u(x,-y)$ in $\Omega \setminus \overline {\left(D_{L 0} \cup D_{R 0}\right) }$. It is worth noting that if $u_{\alpha}$ and $u_{\beta}$ are the solutions for the same harmonic function $H$, then there is a constant $c$ such that  $ u_{\alpha}  = u_{\beta}  + c $ and $$ \nabla u_{\alpha}  = \nabla u_{\beta} $$ in $ \mathbb{R}^2 \setminus \overline {\cup_{n=-\infty}^{\infty}\left(D_{L n} \cup D_{R n}\right)}$.

 \par In this paper, we establish some asymptotes and  optimal blow-up rates for $\nabla u $ in two kinds of  narrow regions in between $D_{Ln}$ and $D_{Rn}$, and in between $D_{Rn}$ and $D_{Rn+1}$. Theorem \ref{upp_thm} provides asymptotes with a  coefficient and an upper bound of the coefficient, and moreover Theorem \ref{low_thm} presents a specific asymptote with a lower bound for the coefficient in the case of a linear function $$ H (x,y) = ax +by$$ to get the optimality of the gradient estimates.

\begin{thm}  \label {upp_thm}  For any harmonic function $H$ with \eqref{H_condi}, let $u$ be a solution to \eqref {main} with the condition \eqref{main_condi}. Let $N_v$ be a  narrow vertical region in between $D_{L0}$ and $D_{R0}$, and let $N_h$ be a  narrow horizontal region in between $D_{R0}$ and $D_{R1}$, defined as
\begin{align*}
&N_v= \left\{ (x,y) ~\Big |~ |x| <  1 +\frac {\epsilon} 2-\sqrt {1 - y^2}\mbox { and }|y|<\frac {\sqrt 3} 2  \right\}   \mbox { and}\\
&N_h=\left\{ (x,y) ~\Big |~ \left|y-1-\frac \delta 2 \right| <  1 +\frac {\delta} 2-\sqrt {1 - \left(x-1-\frac \epsilon 2  \right)^2}\mbox { and }\left|x-1-\frac \epsilon 2 \right|<\frac {\sqrt 3} 2  \right\}.
\end{align*}  Then, there exist constants $ \mu$ and $\lambda$ such that  
\begin{align}&  \nabla u (x,y) =  \lambda  \left( \frac {1}{ \delta + \left(x-1-\frac {\epsilon}2 \right)^2} \left( 0,1\right) +  \frac 1 {\sqrt \delta} S (x,y) (1,0) \right) + R_1 (x,y)~&&\mbox{for any } (x,y)\in N_h, \label{upp_thm_1st_inequ}\\
&  \nabla u (x,y) =  \mu  \frac {\sqrt \epsilon}{ \epsilon + y^2 }   \left(1,0\right) + R_2 (x,y) ~&&\mbox{for any } (x,y)\in N_v. \label{upp_thm_2nd_inequ}\end{align} Here, the constants $\lambda$ and $\mu$ satisfy 
$$ \lambda = H\left(0,1  \right) - H\left(0,-1 \right) $$ and $$ |\mu |\leq C \norm {H}_{L^\infty (\Omega_4)}, $$ the function $S$ is defined as  $$ S(x,y)  = -2 \sqrt {\delta}~  \frac { \left(x-1-\frac \epsilon 2 \right)\left(y-1-\frac \delta 2 \right)}{\left(\left(x-1-\frac \epsilon 2 \right)^2 + \delta\right)^2} $$ with $$ \norm { S (x,y)  }_{L^{\infty} (N_h)} \leq 2 ~\mbox{and}~ S\left(1+\frac \epsilon 2, 1+\frac \delta 2 \right) = 0 ,$$ and the remainder terms  $R_1$ and $R_2$ are bounded as 
 $$\norm {R_1}_{L^\infty (N_h)}  + \norm {R_2}_{L^\infty (N_v)} \leq C \norm {H}_{L^\infty (\Omega_4)} $$  for a constant $C$ regardless of $\epsilon$ and $\delta$.
\end{thm}  It is worth noting that the constant $\mu$ depends on $\delta$ and $\epsilon$ in this paper, even though it is bounded regardless of $\delta$ and $\epsilon$. The proof of Theorem \ref {upp_thm} is given in Section \ref{upp_thm_sec}, based on the results in Section \ref{EPD_sec}.

\begin{rem}  We consider the behavior of $\nabla u$ in the domain $$ \left (-1-\frac {\epsilon} 2 + \frac {\sqrt 3 } 2,1+\frac {\epsilon} 2  -\frac {\sqrt 3 } 2 \right)\times\left( \frac {\sqrt 3 } 2,2+\delta-  \frac {\sqrt 3 } 2 \right) $$ which doesn't belong to  $N_v$ and $N_h$.  Theorem \ref {upp_thm} provides the boundedness of $|\nabla u|$ on its rectangular boundary regardless of $\epsilon$ and $\delta$.  By the maximum principle,  $|\nabla u|$ is bounded in the domian regardless of $\epsilon$ and $\delta$. Combined with Theorem \ref {upp_thm} again, an asymptote for $\nabla u$ in $$\left(\left(-\frac {\sqrt 3} 2 ,\frac {\sqrt 3} 2\right) \times \mathbb{R} \right)\setminus \overline {\cup_{n=-\infty}^{\infty}\left(D_{L n} \cup D_{R n}\right)}$$ is also obtained,  since the gradient $\nabla u$ is periodic with period $2+\delta$ in the $y$ direction. 
\end{rem}

\begin{thm} \label {low_thm} Let $N_v$, $N_h$ and $S$ be as given in Theorem \ref{upp_thm}. Assume that  $H$ is a linear function given as $$H(x,y) = a x + b y $$ for any $(x,y) \in \mathbb{R}^2$ and  $u$ is a solution to \eqref {main} with the condition \eqref{main_condi} for $H$.  Then,
\beq  \nabla u (x,y) =  2 b   \left( \frac {1}{ \delta + \left(x-1-\frac {\epsilon}2 \right)^2} \left( 0,1\right) +   \frac 1 {\sqrt \delta} S(x,y) (1,0) \right) + R_1 (x,y)~\mbox{for any } (x,y)\in N_h,  \label {low_thm_1st_eq}\eeq and there is a constant $\mu_0$ such that \beq \nabla u (x,y) = a  \mu_0 \frac {\sqrt \epsilon}{ \epsilon + y^2 }   \left(1,0\right) + R_2 (x,y) \mbox{ for any } (x,y)\in N_v \label {low_thm_2nd_eq}.\eeq  Here, the constant $\mu_0$ satisfies \beq \frac 1 C < \mu_0  < C \label {low_thm_3rd_eq}\eeq and  the remainder terms $R_1$ and $R_2$ are bounded as 
 $$\norm {R_1}_{L^\infty (N_h)}  + \norm {R_2}_{L^\infty (N_v)} \leq C (|a|+ |b|)$$ for a constant $C$ regardless of $\epsilon$ and $\delta$. \end{thm} The proof of Theorem \ref {low_thm} is presented in Section \ref{low_thm_sec}.

 \par  From now on,  the symbols $C$ and $C_n$ denote the constants regardless of  small $\epsilon >0$ and $\delta >0$ for $n=1,2,\cdots$.

\begin{rem} \label {just_remark} Theorems \ref {upp_thm} and \ref {low_thm} provide the generic blow-up rates of $\nabla u $ in between $ D_{Ln}$ and $ D_{Rn}$, and in between $D_{Ln}$ and $D_{Ln+1}$ which are  $$\frac 1 {\sqrt \epsilon} \mbox { and }\frac 1 {\delta},$$ respectively.

\end{rem}

\begin{cor} \label {cor_excercise} For a  harmonic function $H$ with a periodic gradient as above, let $w$ be a solution to the equation 
$$ \Delta w = 0 \mbox{ in }  \mathbb{R}^2 \setminus \overline {\cup_{n=-\infty}^{\infty} D_{R n} } $$ with the conditions \beq
 \begin{cases}
 \ds  w = c    &\mbox{ on } \p D_{R 0}  ,\\ \nm  \ds \int_{\p D_{R 0}} \p_{\nu} w ds =0,&\\
  \nm \ds  \int_{\Omega \setminus \overline { D_{R0}}} |\nabla (w-H)|^2 dxdy < \infty,  & \\  \nm \ds  \nabla w(x,y) =  \nabla w(x,y + 2 + \delta) & \mbox { for any } (x,y)\in \mathbb{R}^2 \setminus \overline {\cup_{n=-\infty}^{\infty}D_{R n}} ,\end{cases}\notag
 \eeq where $c$ is a constant depending on $H$, $\epsilon$. Then, there exist a constant $\lambda$ such that  
$$ \nabla  w (x,y) =  \lambda    \left( \frac {1}{ \delta + \left(x-1-\frac {\epsilon}2 \right)^2} \left( 0,1\right) + \frac 1 {\sqrt {\delta} }  S (x,y) (1,0) \right) + R(x,y)~~~~\mbox{for any } (x,y)\in N_h, $$ where the constant $ \lambda = H\left(0,1  \right) - H\left(0,-1 \right) $ and the remainder term $R$ is bounded as 
 $$\norm {R}_{L^\infty (N_h)} \leq C \norm {H}_{L^\infty (\Omega_4)} $$ for a constant $C$ regardless of $\delta >0$. 

\end{cor}

This paper is organized as follows: Section \ref {EPD_sec} provides the potential differences of $u$ between $\p D_{R n}$  and $\p D_{L n}$, and between $\p D_{R n}$  and $\p D_{R n+1}$. In Section \ref{upp_thm_sec}, two asymptotes \eqref{upp_thm_1st_inequ} and \eqref {upp_thm_2nd_inequ} for $\nabla u$ in Theorem  \ref {upp_thm} result from  the potential differences. In Section \ref{low_thm_sec}, we establish more descriptive asymptotes of  $\nabla u$  for $H(x,y)=ax + by$ to  prove Theorem \ref {low_thm} and to  get the optimality of the blow-up rates in Remark \ref{just_remark}.  The proof of Corollary \ref {cor_excercise}  is left as an exercise for the reader, since it is much the same as the proof of \eqref{upp_thm_1st_inequ} and Proposition \ref{prop_poten_differ_2}.

\section{Estimates for potential differences} \label{EPD_sec}
The potential differences play  important roles in establishing the asymptotes and estimates for the gradient $\nabla u$.  Once the potential difference is estimated, the methods in \cite{KLY,KLeeY} are modified to obtain the asymptotes. This section thus provides the estimates for potential differences $u |_{D_{R n+1}} - u |_{D_{Rn}}$, $u |_{D_{L n+1}} - u |_{D_{Ln}} $ and $u |_{D_{R n}} - u |_{D_{Ln}} $.

\begin{prop}\label{prop_poten_differ_2}  Let $u$ be a solution to \eqref {main} with the condition \eqref{main_condi} for any harmonic function $H$ with a periodic gradient satisfying \eqref{H_condi}.  Then, the potential differences  between $\p D_{L n}$  and $\p D_{L n+1}$, and between $\p D_{R n}$  and $\p D_{R n+1}$ are obtained as 
$$ u |_{D_{R n+1}} - u |_{D_{Rn}} =    u |_{D_{L n+1}} - u |_{D_{Ln}} =  H\left(0,1 + \frac \delta 2 \right) - H\left(0,-1 - \frac \delta 2 \right)$$  for any $n =0, \pm 1, \pm 2, \cdots$.
\end{prop}

\pf We begin by proving that \beq (u - H) (x,y) - (u - H) \left(x,y -2-\delta \right) = 0 \label{proof_first_prop}\eeq for any $(x,y) \in \mathbb{R}^2 \setminus \overline {\cup_{n=-\infty}^{\infty}\left(D_{L n} \cup D_{R n}\right)} $.  The gradient $\nabla (u - H)$ is a periodic function with period $2+\delta$ in the $y$ direction as given in \eqref{H_condi} and \eqref{main_condi}. Then, $\nabla (u - H) (x,y) - \nabla (u - H) \left(x,y -2-\delta \right) = (0,0) $ for any $(x,y) \in \mathbb{R}^2 \setminus \overline {\cup_{n=-\infty}^{\infty}\left(D_{L n} \cup D_{R n}\right)} $. There exists a constant $d>0$ such that $$ (u - H) (x,y) - (u - H) \left(x,y -2-\delta \right) = d $$ for any $(x,y) \in \mathbb{R}^2 \setminus \overline {\cup_{n=-\infty}^{\infty}\left(D_{L n} \cup D_{R n}\right)} $. By the Jensen's inequality, every $x \in [3,\infty)$ has the upper bound for $|d|^2$ as  $$|d |^2 \leq  \left(\int_{-1-\frac {\delta} 2 } ^{1+\frac {\delta} 2 } | \p_y (u-H) (x,y)| dy \right)^2  \leq ( 2 +\delta)\int_{-1-\frac {\delta} 2 } ^{1+\frac {\delta} 2 } | \nabla (u-H) (x,y)|^2 dy .$$
Since $\int_{\Omega \setminus \overline { D_{L0} \cup D_{R0}}} |\nabla (u-H)|^2 dxdy < \infty$, it follows from the Fubini's theorem that $$d = 0$$ implying  \eqref{proof_first_prop}.

\par The periodic property \eqref{H_condi} implies that \beq H (x,y) - H \left(x,y -2-\delta \right) = H\left(0,1 + \frac \delta 2 \right) - H\left(0,-1 - \frac \delta 2 \right)\label{H=H-H}\eeq  for any $(x,y) \in \mathbb{R}^2 $, since $\nabla \left(H (x,y) - H \left(x,y -2-\delta \right) \right) = (0,0) $. The equality \eqref{proof_first_prop} yields this proposition as follows: 
\begin{align*}
 &H\left(0,1 + \frac \delta 2 \right) - H\left(0,-1 - \frac \delta 2 \right) = H (x,y) - H \left(x,y -2-\delta \right)  \\&= H (x,y) - H \left(x,y -2-\delta \right) + \left( (u - H) (x,y) - (u - H) \left(x,y -2-\delta \right) \right)\\&=   u  (x,y) - u  \left(x,y -2-\delta \right) 
\end{align*} for any $(x,y) \in \mathbb{R}^2 \setminus \overline {\cup_{n=-\infty}^{\infty}\left(D_{L n} \cup D_{R n}\right)} $.   This implies that $$H\left(0,1 + \frac \delta 2 \right) - H\left(0,-1 - \frac \delta 2 \right) = u |_{D_{L n+1}} - u |_{D_{Ln}} = u |_{D_{R n+1}} - u |_{D_{Rn}}. $$
\qed

\begin{prop}\label{prop_poten_differ_1} Let $u$ be the solution to \eqref {main}   satisfying \eqref{main_condi} for   $$H(x,y)=x$$ for any $(x,y) \in \mathbb{R}^2$. Then, the potential difference between $\p D_{L n} $ and $\p D_{R n} $ is estimated as 
$$ \frac 1 C \sqrt \epsilon <  u |_{D_{R n}} - u |_{D_{L n}} < C \sqrt \epsilon,$$ and there are no potential differences between $\p D_{L n}$  and $\p D_{L n+1}$, and between $\p D_{R n}$  and $\p D_{R n+1}$, i.e.,$$   u |_{D_{L n+1}} - u |_{D_{L n}}= u |_{D_{R n+1}} - u |_{R_{L n}}= 0$$ for any $n =0, \pm 1, \pm 2, \cdots$.
\end{prop}
The proof of the proposition is presented in Subsection \ref{subsec_1st_prop}.

The potential difference of $u$ between $\p D_{L0}$ and $\p D_{R0}$ can be expressed as an integral containing $\p_{\nu} \phi$ in Proposition \ref{phi_lemma} motivated by the method in \cite{Y, Y2}.  Following lemma is used to modify the idea for the proposition.

\begin{lem} \label{lem2-4-h} Let $h$ be a harmonic function as $$
 \begin{cases}
 \ds\triangle h = 0   &\mbox{ in }(4,\infty)\times\left(1-\frac \delta 2 , 1+\frac \delta 2\right),  \\ \nm  \ds \p_y h = 0     & \mbox{on }(4,\infty)\times\left \{1-\frac \delta 2 , 1+\frac \delta 2\right\} ,   \\   \nm \ds  \int_{(4,\infty)\times\left(1-\frac \delta 2 , 1+\frac \delta 2\right)}  |\nabla h|^2 dxdy & < \infty.  \end{cases}
$$ Then, $$ \sup_{y \in \left(1-\frac \delta 2 , 1+\frac \delta 2\right)} | h (x,y) |= O(1) ~\mbox{and}~\sup_{y \in \left(1-\frac \delta 2 , 1+\frac \delta 2\right)}  \left|\nabla h (x,y) \right|=O\left(\exp\left(-{\frac {\pi}{2+\delta}} x \right)\right)$$  as $ x\rightarrow \infty $.\end{lem}
\pf The function $h$ can be express as $$h(x,y) = \sum_{n=0}^{\infty} a_n \cos\left({\frac {n\pi}{2+\delta}} \left(y+1+\frac \delta 2 \right)\right)\exp\left(-{\frac {n\pi}{2+\delta}} x \right).$$ The estimates can be obtained immediately. \qed

\begin{prop} \label{phi_lemma} There exists the harmonic function $\phi$ defined in  $\Omega \setminus \overline {\left(D_{L0} \cup D_{R0}\right)} $ with the following conditions as
\beq
 \begin{cases}
 \ds \partial_{\nu} \phi= 0   &\mbox{on } \mathbb{R} \times \left\{y~\Big|~ y= -1 - \frac {\delta}  2 ~\mbox{or}~y= 1 + \frac {\delta}  2\right\},  \\ \nm  \ds \phi =  c_0    & \mbox{on } \partial  D_{R0},   \\   \nm \ds  \phi =- c_0  & \mbox{on } \partial  D_{L0}, \\ \nm \ds \int_{\partial D_{R0}} \partial_{\nu} \phi  ds =  -\int_{\partial D_{L0}} & \partial_{\nu}  \phi  ds =  1 , \\  \nm \ds  \int_{\Omega \setminus \overline {D_{L0} \cup D_{R0}}} \left|\nabla \phi \right |^2 dxdy & < \infty, \end{cases}
 \eeq where $c_0$ is a proper constant depending on $\epsilon$ and $\delta$.  If $u$ is a solution to \eqref {main} satisfying  the condition \eqref{main_condi} for  any harmonic function $H$ with \eqref{H_condi}, then \beq u \Big|_{\partial D_{R0}} - u \Big|_{\partial D_{L0}}  = \int_{\partial D_{L0} \cup \partial D_{R0}} H \partial_{\nu} \phi  ds . \label{pro_equal} \eeq
\end{prop}

\pf First, we prove the existence of $\phi$. By the Lax-Milgram theorem, there exists the harmonic function $\varphi_0$ defined in $\Omega \setminus \overline {\left(D_{L0} \cup D_{R0}\right)} $ with conditions:$$
 \begin{cases}
 \ds \partial_{\nu} \varphi_0= 0   &\mbox{on } \mathbb{R} \times \left\{y~\Big|~ y= -1 - \frac {\delta}  2 ~\mbox{and}~y= 1 + \frac {\delta}  2\right\},  \\ \nm  \ds \varphi_0 = 1    & \mbox{on } \partial  D_{R0},   \\   \nm \ds  \varphi_0 =-1  & \mbox{on } \partial  D_{L0}, \\  \nm \ds  \int_{\Omega \setminus \overline {D_{L0} \cup D_{R0}}} \left|\nabla \varphi_0 \right |^2 dxdy & < \infty. \end{cases}
$$ We can constuct a bijective conformal mapping $\varPhi: B_1 (0,0) \rightarrow \Omega$ such that 
$$ \begin{cases}
 \ds \triangle \varphi_0(\varPhi) = 0   &\mbox{in } B_1 (0,0)\setminus \overline {\varPhi^{-1} (D_{L0})\cup \varPhi^{-1} (D_{R0})} ,  \\ \nm  \ds \varphi_0(\varPhi) = 1    & \mbox{on } \partial  \varPhi^{-1} ( D_{R0}),   \\   \nm \ds  \varphi_0 (\varPhi) =-1  & \mbox{on } \partial \varPhi^{-1} ( D_{L0}), \\  \nm \ds  \varphi_0(\varPhi) \mbox{ belongs to } & H^1 (B_1 (0,0)\setminus \overline {\varPhi^{-1} (D_{L0})\cup \varPhi^{-1} (D_{R0})}) , \end{cases}$$ since $$ \int_{B_1 (0,0)\setminus \overline {\varPhi^{-1} (D_{L0})\cup \varPhi^{-1} (D_{R0})}} \left|\nabla ( \varphi_0 (\varPhi)) \right |^2 d\xi d\eta = \int_{\Omega \setminus \overline {D_{L0} \cup D_{R0}}} \left|\nabla \varphi_0 \right |^2 dx dy < \infty.$$
By the maximum principle,  $\varphi_0$ has the maximal value $1$ on $\partial  D_{R0}$ and also has the minimal value $-1$ on $\partial  D_{L0}$, since $\varphi_0(\varPhi)$ is a harmonic function defined in a bounded domain.  By the Hopf's lemma,  $$  \int_{\partial D_{R0}} \partial_{\nu} \varphi_0  ds =  -\int_{\partial D_{L0}}  \partial_{\nu}  \varphi_0  ds >0  ,  $$ since the normal vector $\nu$ points toward the inside of $D_{L0}$ or  $D_{R0}$. Then, $$ \phi  = \frac 1 { \int_{\partial D_{R0}} \partial_{\nu} \varphi_0  ds } ~\varphi_0 ,$$ which means the existence of $\phi$. In addition, it can be easily shown that \beq  \phi  (x,y) = \phi  (x,-y) \label{varphi_even} \eeq for $(x,y) \in \Omega \setminus \overline {D_{L0} \cup D_{R0}}$.

\par Second, we prove the equality \eqref{pro_equal}. From definition, $u$ is constant on each of boundaries $\partial D_{R0}$ and $\partial D_{L0}$, and $\int_{\partial D_{R0}} \partial_{\nu} \phi  ds =  -\int_{\partial D_{L0}}  \partial_{\nu}  \phi  ds =  1.$ Thus, \begin{align*} u \Big|_{\partial D_{R0}} - u \Big|_{\partial D_{L0}}  & = \int_{\partial D_{L0} \cup \partial D_{R0}}  u \partial_{\nu} \phi  ds \\ 
&= \int_{\partial D_{L0} \cup \partial D_{R0}}  H \partial_{\nu} \phi  ds + \int_{\partial D_{L0} \cup \partial D_{R0}} ( u- H) \partial_{\nu} \phi  ds. \end{align*}
We shall use the divergence theorem to prove that $$ \int_{\partial D_{L0} \cup \partial D_{R0}} ( u- H) \partial_{\nu} \phi = 0.$$ This immediately results in the desirable equality \eqref{pro_equal}. To use the divergence theorem, we define $\tilde u $ and $\tilde H$ as even functions with respect to $y$ as $\tilde u (x,y) = \frac 1 2 \left( u (x,y) + u (x,-y)\right)$ and 
$\tilde H (x,y) = \frac 1 2 \left( H (x,y) + H (x,-y)\right).$  By the periodic property of $\nabla u$,   $\tilde u - \tilde H$ has zero Neumann data on two horizontal boundaries $\p \Omega$ so that   
$$\p_{y} \left( \tilde u - \tilde H \right) \left(x, 1+\frac \epsilon 2\right) =  \p_{y}  \left( \tilde u - \tilde H \right)  \left(x, -1-\frac \epsilon 2\right)  = 0  $$ for any $x \in \mathbb{R}$. From definition of $u$, $$ \int_{\Omega \setminus \overline { D_{L0} \cup D_{R0}}} \left|\nabla \left( \tilde u - \tilde H \right)\right|^2 dxdy < \infty.$$ By Lemma \ref{lem2-4-h}, 
 $$  \tilde u (x,y)- \tilde H(x,y)  = O(1)~\mbox{and}~ \nabla \left(\tilde u (x,y)- \tilde H(x,y) \right) = O\left(\exp -|x|\right)$$  as $|x|\rightarrow \infty $, and $\phi$ and $\nabla \phi$ also show the same behaviors as $|x|\rightarrow \infty $.  Thus, we can use the divergence theorem so that by \eqref {varphi_even}, \begin{align*}&\int_{\partial D_{L0} \cup \partial D_{R0}} ( u- H) \partial_{\nu} \phi  ds \\&= \int_{\partial D_{L0} \cup \partial D_{R0}} ( \tilde u- \tilde H) \partial_{\nu} \phi ds\\ &= \int_{\partial \left(\Omega \setminus \overline { D_{L0} \cup D_{R0}}\right)} ( \tilde u- \tilde H) \partial_{\nu} \phi  ds  \\&= \int_{\partial \left({ D_{L0} \cup D_{R0}}\right)} \phi  \partial_{\nu} ( \tilde u- \tilde H) ds  \\&= \phi |_{\partial D_{L0}}  \int_{\partial D_{L0}}  \partial_{\nu} ( \tilde u- \tilde H) ds  +\phi |_{\partial D_{R0}}  \int_{\partial D_{R0}}  \partial_{\nu} ( \tilde u- \tilde H) ds =0.\end{align*} Thus, we have done it.
\qed

\subsection {Proof of Proposition \ref{prop_poten_differ_1} }\label{subsec_1st_prop}
In this subsection, we suppose that $$ H(x,y)= x $$ for any $(x,y) \in \mathbb{R}^2$.  The function $u$ is the solution to \eqref {main}   satisfying \eqref{main_condi} for $H$.

\par The integral equation \eqref {pro_equal} is mainly used to estimate the potential difference of $u$ between $\partial D_{L0}$ and $\partial D_{R0}$. In \eqref {phi=widephi+v}, the function $\phi$ is construct by a series of $\phi_n$ whose property has been well known, and which is also given explicitly as in \eqref {eqn_phi_n}.

\par First, some maximum principles related to $\phi$ are considered  in Lemmas \ref {lemma_phi_R} and \ref {lemma_rho} before constructing  $\phi$.  Let $\Omega_R$ be the right-hand side of $\Omega$ as
\beq \Omega_R = (0,\infty) \times \left(-1 - \frac {\delta}  2,~ 1 + \frac {\delta}  2 \right). \label{def_Omega_R}\eeq

\begin{lem} \label {lemma_phi_R} There exists a harmonic function $\phi_R$ defined on $\Omega_R \setminus \overline {D_{R0}} $ with the conditions \beq
 \begin{cases}
 \ds \partial_{\nu} \phi_R = 0 \quad & \mbox{on } (0,\infty) \times \left\{y~\Big|~ y= -1 - \frac {\delta}  2 ~\mbox{or}~y= 1 + \frac {\delta}  2\right\},  \\ \nm \ds
 \phi_R = 0   \quad &\mbox{on } \left\{0\right\}  \times \left(-1 - \frac {\delta}  2,~ 1 + \frac {\delta}  2 \right) ,\\ \nm \ds \phi_R = 1   \quad & \mbox{on } \partial  D_{R0},\\  \nm \ds \int_{\Omega_R \setminus \overline {D_{R0}}} \left|\nabla \phi_R \right |^2  dxdy & < \infty. \end{cases} \label{lemma_phi_R_condi}
 \eeq  The function $\phi_R$ has the extreame values only on the boundary so that $$0< \phi_R < 1 ~\mbox{in } \Omega_R \setminus \overline {D_{R0}} .$$
\end{lem} 
\begin{rem} It is obvious that 
$$ \phi_R = \frac 1{c_0} \phi $$ in $\Omega_R \setminus \overline {D_{R0}} $, where the constant $c_0$ is given in Proposition \ref {phi_lemma}.

\end{rem}

\pf By the Lax-Milgram theorem, there exists the unique harmonic function $\phi_R$ defined on $\Omega_R \setminus \overline {D_{R0}} $ with the boundary condition \eqref{lemma_phi_R_condi}. As mentioned in the remark above, the existence of $\phi_R$ results immediately from $\phi$ given  in Proposition \ref{phi_lemma}  due to $\phi_R = \frac 1{c_0} \phi$.

\par We use a conformal map to prove that $0< \phi_R < 1$ in $\Omega_R \setminus \overline {D_{R0}} .$ Let $$B_1 ^+ (0,0) = \left\{ (\xi,\eta) ~|~\xi^2 + \eta^2 < 1 \mbox{ and }\xi >0 \right\}.$$ There exists  a bijective conformal mapping $\varPhi_R: B_1^+ (0,0) \rightarrow \Omega_R$ such that 
$$ \begin{cases}
 \ds \triangle \phi_R(\varPhi_R) = 0   &\mbox{in } B_1 ^+ (0,0)\setminus \overline { \varPhi_R ^{-1} (D_{R0})} ,  \\ \nm  \ds \phi_R(\varPhi_R) = 1    & \mbox{on } \partial  \varPhi_R^{-1} ( D_{R0}),   \\   \nm \ds  \phi_R (\varPhi_R) = 0  & \mbox{on }  \xi^2 + \eta^2 = 1 \mbox{ and  } \xi >0 ,\\   \nm \ds \partial_{\nu} \left(\phi_R (\varPhi_R)\right)(0,\eta) = 0  & \mbox{on }  |\eta|<1 , \\  \nm \ds  \phi_R(\varPhi_R) \mbox{ belongs to } & H^1 (B_1 (0,0)\setminus \overline { \varPhi_R ^{-1} (D_{R0})}) . \end{cases}$$  By the maximal principle, $0 < \phi_R(\varPhi_R) <1 $ in the bounded domain $B_1 ^+ (0,0)\setminus \overline { \varPhi_R ^{-1} (D_{R0})}$. This implies that   $$0< \phi_R < 1 ~\mbox{in } \Omega_R \setminus \overline {D_{R0}} .$$\qed

\par The following lemma is derived easily by an argument, analogous to Lemma \ref{lemma_phi_R}, where  the same function $\varPhi_R: B_1^+ (0,0) \rightarrow \Omega_R$ is used.
 
\begin{lem} \label{lemma_rho} Let $\rho$ be a harmonic function defined in $\Omega_R \setminus \overline {D_{R0}} $ with the boundary conditions: \beq
 \begin{cases}
 \ds \partial_{\nu} \rho = 0 \quad & \mbox{on } (0,\infty) \times \left\{y~\Big|~ y= -1 - \frac {\delta}  2 ~\mbox{or}~y= 1 + \frac {\delta}  2\right\},  \\ \nm \ds
 \rho = 0   \quad &\mbox{on } \left\{x~\Big|~ x=0\right\}  \times \left(-1 - \frac {\delta}  2,~ 1 + \frac {\delta}  2 \right) ,\\ \nm \ds \rho >0  \quad & \mbox{on } \partial  D_{R0}, \end{cases}
 \eeq and also satisfying 
$$ \int_{\Omega_R \setminus \overline {D_{R0}}} \left|\nabla \rho \right |^2  dxdy < \infty. $$ Then, $$0< \rho ~\mbox{in } \Omega_R \setminus \overline {D_{R0}} .$$
\end{lem}

\vskip10pt
\par Second, we use a series of functions $\phi_n$ to express $\phi$, given in Proposition \ref{phi_lemma}. Here, the harmonic function $\phi_n $ satisfies \beq
 \begin{cases}
 \ds \Delta \phi_n  = 0 \quad  & \mbox{in } \mathbb{R}^2 \setminus \overline{D_{Ln} \cup D_{Rn}}\\
 \nm \ds  \phi_n = \mbox{ a constant}   \quad & \mbox{on } \partial  D_{Ln}, \\ \nm  \ds \phi_n = - \phi_n \Big|_{ \partial  D_{Ln}}   ~ \quad & \mbox{on } \partial  D_{Rn},   \\ \nm \ds \int_{\partial D_{Rn}} \partial_{\nu} \phi_n  ds =  -\int_{\partial D_{Ln}} \partial_{\nu}  \phi_n  ds &=  1 , \\  \nm \ds   \phi_n ({\bf x}) = O\left({\frac 1 {|{\bf x}|}}\right) & \mbox{as } {\bf x} \rightarrow \infty \end{cases}\label{eqn_phi_n_before} \eeq  for any integer $n$. Then, the function $\phi_n$ is expressed explicitly as
\beq \phi_n = \frac 1 {2\pi}\big(\log \left| {\bf x} - (-p, n(2+\delta)) \right| - \log \left| {\bf x} - (p, n(2+\delta)) \right|\big) \label{eqn_phi_n} \eeq where $$p = \sqrt {\epsilon} + O(\epsilon) $$ for small $\epsilon>0$. Refer to \cite{KLY,LY} for details. We define the sum $\widetilde {\phi}$ of the series  of $\phi_n$  in the manner as  $$ \widetilde {\phi} = \lim_{N \rightarrow \infty}\sum_{n= -N } ^{N} \phi_n.$$  The function $\widetilde {\phi} $ is well defined in $\Omega \setminus \overline {D_{L0} \cup D_{R0}}$ by the help of  a neutralization reaction between $\phi_n $ and $ \phi_{-n}$, and satisfies \beq
 \begin{cases}
 \ds \triangle \widetilde { \phi} = 0   &\mbox{in } \Omega \setminus \overline {D_{L0} \cup D_{R0}}, \\ \nm \ds \partial_{\nu}\widetilde { \phi}= 0   &\mbox{on } \p \Omega =  \mathbb{R} \times \left\{y~\Big|~ y= -1 - \frac {\delta}  2 ~\mbox{or}~y= 1 + \frac {\delta}  2\right\},  \\ \nm \ds \int_{\partial D_{R0}} \partial_{\nu} \widetilde {\phi}  ds =  -\int_{\partial D_{L0}} & \partial_{\nu} \widetilde { \phi}  ds =  1 , \\  \nm \ds  \int_{\Omega \setminus \overline {D_{L0} \cup D_{R0}}} \left|\nabla \widetilde {\phi} \right |^2 dxdy & < \infty. \end{cases} \notag
 \eeq Here, $\widetilde { \phi} $ is not constant on each of $\p D_{L0}$ and $\p D_{R0}$, and $\Omega$ is as given at \eqref{def_Omega}. There exists a harmonic function $\widetilde {v}$ defined in $\Omega \setminus \overline {\left(D_{L0} \cup D_{R0}\right)} $ with conditions:
\beq
 \begin{cases}
 \ds \partial_{\nu}\widetilde {v} = 0   &\mbox{on } \p \Omega =  \mathbb{R} \times \left\{y~\Big|~ y= -1 - \frac {\delta}  2 ~\mbox{or}~y= 1 + \frac {\delta}  2\right\},  \\ \nm  \ds \widetilde {\phi}+ \widetilde {v} =  \mbox{ a constant }  \widetilde {c}    & \mbox{on } \partial  D_{R0},   \\   \nm \ds  \widetilde {\phi} + \widetilde {v} =- \widetilde {c}  & \mbox{on } \partial  D_{L0}, \\ \nm \ds \int_{\partial D_{R0}} \partial_{\nu} \widetilde {v}  ds  =  \int_{\partial D_{L0}} & \partial_{\nu} \widetilde {v}  ds =  0, \\  \nm \ds  \int_{\Omega \setminus \overline {D_{L0} \cup D_{R0}}} \left|\nabla \widetilde {v} \right |^2 dxdy & < \infty \end{cases}
 \eeq for a proper constant $\tilde c $. The existence of $\tilde v$ is derived in the same way as $u$ to \eqref {main} for a given $H$ in the introduction. In another way, the existence of $\tilde v$ is also derived  from  the existence of $\phi$ shown in Proposition \ref{phi_lemma}, since  $ \widetilde {\phi} + \widetilde {v}$ satisfies all conditions of  $\phi$ so that  \beq \phi = \widetilde {\phi} + \widetilde {v}  \label{phi=widephi+v}.\eeq  
 
The function $ \phi$ can be decomposed into three functions as 
 $$ \phi = \alpha \phi + \left(\widetilde {\phi} -  \alpha \phi \right) + \widetilde {v}, $$
where the positive constant $\alpha$ is defined as 
$$ \alpha = \frac {\widetilde {\phi}{\left(\frac \epsilon 2 ,0 \right)}- \widetilde {\phi}{\left(-\frac \epsilon 2 ,0 \right)} } { {\phi} \Big|_{\p D_{R0}}-  {\phi}\Big|_{\p D_{L0}} }.$$

\vskip 20pt
\begin{lem} \label{last_lemma_ok} There is a constant $C$ such that 
$$0 < 1- \alpha \leq C \sqrt{\epsilon}$$ for small $\epsilon >0 $.
\end{lem}
\pf First, we show that \begin{align} 0 &<\left(\widetilde {\phi} -  \alpha \phi + \widetilde {v} \right) \Big|_{\p D_{R0}} - \left(\widetilde {\phi} -  \alpha \phi + \widetilde {v}\right) \Big|_{\p D_{L0}} \notag \\&= \int_{\p D_{L0} \cup \p D_{R0}} \left(\widetilde {\phi} -  \alpha \phi  \right) \partial_{\nu} \phi ~ds \notag \\&\leq  C_1 \sqrt {\epsilon} \int_{\p D_{L0} \cup \p D_{R0}}  x \partial_{\nu} \phi ~ds . \label{1steq_last_lemma} \end{align} for a positive constant $C_1$. In the same way as Proposition \ref{phi_lemma}, the integation by parts yields \begin{align*}&\int_{\p D_{L0} \cup \p D_{R0}}  \tilde v \p_{\nu} \phi  ds = \int_{\p D_{L0} \cup \p D_{R0}}  \phi  \p_{\nu} \tilde v   ds \\&= \left(\phi \Big|_{\p D_{L0}}\right)\int_{\p D_{L0}}   \p_{\nu} \tilde v   ds + \left( \phi \Big|_{\p D_{R0}}\right) \int_{ \p D_{R0}}   \p_{\nu} \tilde v   ds  = 0. \end{align*}  We thus have the equality \begin{align*} & \left(\widetilde {\phi} -  \alpha \phi + \widetilde {v} \right) \Big|_{\p D_{R0}} - \left(\widetilde {\phi} -  \alpha \phi + \widetilde {v}\right) \Big|_{\p D_{L0}} =           \int_{\p D_{L0} \cup \p D_{R0}} \left(\widetilde {\phi} -  \alpha \phi + \tilde v \right) \partial_{\nu} \phi ~ds \\ &=   \int_{\p D_{L0} \cup \p D_{R0}} \left(\widetilde {\phi} -  \alpha \phi  \right) \partial_{\nu} \phi ~ds. \end{align*} Meanwhile, calculationg $\widetilde {\phi} (x,y)$ directly from \eqref {eqn_phi_n}, there exists a positve constant $C_2 $ regardless of $\epsilon$ and $\delta$ such that  
\begin{align*}
0 < \widetilde {\phi} (x,y) - \alpha \phi (x,y) \leq C_2  \sqrt \epsilon x ~ \quad \mbox{for }(x,y) \in \p D_{R0} \setminus \{\left(\epsilon / 2 , 0 \right)\} , \\
 C_2\sqrt \epsilon x  \leq  \widetilde {\phi} (x,y)- \alpha \phi (x,y) < 0 ~ \quad \mbox{for }(x,y) \in \p D_{L0}  \setminus \{\left(- \epsilon /2 , 0 \right)\}, 
\end{align*} and $\widetilde {\phi} \left(\frac \epsilon 2 ,0 \right) - \alpha \phi \left(\frac \epsilon 2 ,0 \right) = \widetilde {\phi} \left(-\frac \epsilon 2 ,0 \right) - \alpha \phi \left(-\frac \epsilon 2 ,0 \right)  = 0$. Thus,  $$ 0 <   \int_{\p D_{L0} \cup \p D_{R0}} \left(\widetilde {\phi} -  \alpha \phi  \right) \partial_{\nu} \phi ~ds   \leq  C_2 \sqrt {\epsilon} \int_{\p D_{L0} \cup \p D_{R0}}  x \partial_{\nu} \phi ~ds,$$  since $\partial_{\nu} \phi >0 $ on $\p D_{R0}$ and  $\partial_{\nu} \phi <0 $ on $\p D_{L0}$ by the Hopf lemma. This implies \eqref{1steq_last_lemma}.

\par  Second, the positivity of $1-\alpha$ is derived simply from \eqref{1steq_last_lemma}. We note that $$ \left(\widetilde {\phi} -  \alpha \phi + \widetilde {v} \right) \Big|_{\p D_{R0}} = - \left(\widetilde {\phi} -  \alpha \phi + \widetilde {v} \right) \Big|_{\p D_{L0}}  = (1-\alpha ) \phi \Big|_{\p D_{R0}}$$ and $ \phi \Big|_{\p D_{R0}} >0$. It follows immediately from \eqref{1steq_last_lemma} that  $$1-\alpha > 0.$$

\par Third,  we consider $ \int_{\p D_{L0} \cup \p D_{R0}}  x \partial_{\nu} \phi ~ds$ to estimate $1-\alpha$.  The inequality \eqref{1steq_last_lemma} implies $$ \int_{\p D_{L0} \cup \p D_{R0}}  x \partial_{\nu} \phi ~ds >0.$$ The integral can be decomposed into two terms as 
 \begin{align}
\int_{\p D_{L0} \cup \p D_{R0}}  x \partial_{\nu} \phi ~ds = \alpha \int_{\p D_{L0} \cup \p D_{R0}}  x \partial_{\nu} \phi ~ds + \int_{\p D_{L0} \cup \p D_{R0}}  x \partial_{\nu}  \left(\widetilde {\phi} -  \alpha \phi + \widetilde {v}\right) ~ds.  \label{decomp_xpphi}
\end{align}
By Lemma \ref{lem2-4-h}, the boundedness of $ \int_{\Omega \setminus \overline {D_{L0} \cup D_{R0}}} \left|\nabla \left(\widetilde {\phi} -  \alpha \phi + \widetilde {v} \right) \right |^2 dxdy $ implies the existence of  a constant  $c_1$  such that  $\left(\widetilde {\phi} -  \alpha \phi + \widetilde {v}\right) (x,y)$ converges to a constant $c_1$, or  $-c_1$, as $x$ approaches $\infty$ or $-\infty, $ respectively, and also show that $ \p_{x} \left(\widetilde {\phi} -  \alpha \phi + \widetilde {v}\right) (x,y)$ shrinks exponentially fast, to $0$, as $|x|$ approaches $\infty$.
The integration by parts yields  \begin{align*}
&\left| \int_{\p D_{L0} \cup \p D_{R0}}  x \partial_{\nu}  \left(\widetilde {\phi} -  \alpha \phi + \widetilde {v}\right) ~ds  \right|\\ &= \left| \lim_{l \rightarrow \infty} \int_{\p D_{L0} \cup \p D_{R0} \cup \left(\{x| x = \pm l \} \times (-1 - {\delta}/  2,~1 + {\delta}/  2)\right)} \partial_{\nu}  x   \left(\widetilde {\phi} -  \alpha \phi + \widetilde {v}\right) ~ds  \right| \\ &= \lim_{l \rightarrow \infty}  \left| \int_{ \{ l \} \times (-1 - {\delta}/  2,~1 + {\delta}/  2)}   \widetilde {\phi} -  \alpha \phi + \widetilde {v}~ds   -  \int_{ \{ - l \} \times (-1 - {\delta}/  2,~1 + {\delta}/  2)}   \widetilde {\phi} -  \alpha \phi + \widetilde {v}~ds \right|,
\end{align*} since $  \widetilde {\phi} -  \alpha \phi + \widetilde {v} $ is constant on $\p D_{L0}$ and on $\p D_{R0}$, respectively. Note that $\widetilde {\phi} -  \alpha \phi + \widetilde {v} = (1-\alpha) \phi $,  $(1-\alpha) >0$,  $\phi(x,y) = -\phi(-x,y)$ and $\phi = \frac 1 {c_0} \phi_{R}  $ for $(x,y) \in \Omega_R \setminus \overline {D_{R0}}$, where $\phi_R $ and $c_0>0$ are defined in Lemma \ref{lemma_phi_R} and Proposition \ref{phi_lemma}, respectively. The maximum principle in Lemmas  \ref{lemma_phi_R} and \eqref {1steq_last_lemma} yield \begin{align*} &\left| \int_{\p D_{L0} \cup \p D_{R0}}  x \partial_{\nu}  \left(\widetilde {\phi} -  \alpha \phi + \widetilde {v}\right) ~ds  \right|\\&= 2 \lim_{l \rightarrow \infty}  \left| \int_{ \{ l \} \times (-1 - {\delta}/  2,~1 + {\delta}/  2)}   \widetilde {\phi} -  \alpha \phi + \widetilde {v}~ds  \right|
 \\& \leq (2+\delta) \left( \left(\widetilde {\phi} -  \alpha \phi + \widetilde {v}\right) \Big|_{\p D_{R0}} - \left(\widetilde {\phi} -  \alpha \phi + \widetilde {v}\right) \Big|_{\p D_{L0}}\right) \\ &\leq 3  \int_{\p D_{L0} \cup \p D_{R0}}   \left(\widetilde {\phi} -  \alpha \phi \right)   \partial_{\nu}\phi ~ds  \\ &\leq 3 C_2\sqrt {\epsilon} \int x \partial_{\nu} \phi ds,  \end{align*} since $\delta<1$.
Applying this bound to the decompostion \eqref{decomp_xpphi}, 
$$0<(1-\alpha) \int_{\p D_{L0} \cup \p D_{R0}}  x \partial_{\nu} \phi ~ds \leq 3 C_2 \sqrt{\epsilon} \int_{\p D_{L0} \cup \p D_{R0}}  x \partial_{\nu} \phi ~ds.$$ Since $\int_{\p D_{L0} \cup \p D_{R0}}  x \partial_{\nu} \phi ~ds>0$, we are done.
\qed

\vskip 20pt 

Now, we take the last step to prove Proposition \ref{prop_poten_differ_1}. Calculationg $\widetilde {\phi} (x,y)$ directly from \eqref {eqn_phi_n}, there exists a positve constant $C_* $ regardless of $\epsilon$ and $\delta$ such that
$$ \frac 1 {C_*} \sqrt \epsilon \leq \tilde {\phi} (x,y) \leq C_* \sqrt \epsilon $$ for any $(x,y) \in \p D_{R0}$ containing $\left( \frac {\epsilon}{2},0\right)$. The definition and symmetric property imply $ \alpha = \frac {\tilde \phi \left( \frac {\epsilon}{2},0\right) }{\phi  |_{\p D_{R0}}},$ and Lemma \ref{last_lemma_ok} yields $\frac 1 2 \leq \alpha \leq 2.$ Thus, there is a constant $C_{**}>0$ regardless of $\epsilon$ and $\delta$ such that $$ \frac 1 {C_{**}} \widetilde {\phi} \leq \phi \leq  C_{**} \widetilde {\phi}~\quad  \mbox{on}~ \partial D_{R0}.$$ By Lemma \ref{lemma_rho}, this inequality on the bounday can be extended into  $\Omega_R \setminus \overline{ D_{R0}}$ so that
$$ \frac 1 {C_{**}} \widetilde {\phi} \leq \phi \leq  C_{**} \widetilde {\phi} ~\quad  \mbox{in}~\Omega_R \setminus \overline{ D_{R0}}.$$ By the divergence theorem, 
$$ \int_{\p D_{R0} }  x \p_{\nu} \phi ds = \lim_{l \rightarrow \infty} \int_{-1 - \frac \delta 2}^{1 + \frac \delta 2} \p_{\nu} x \phi (l,y) dy ~=~ \lim_{l \rightarrow \infty} \int_{-1 - \frac \delta 2}^{1 + \frac \delta 2}  \phi  (l,y) dy$$ whose value is  intermediate between $$\frac 1 {C_{**}} \lim_{l \rightarrow \infty} \int_{-1 - \frac \delta 2}^{1 + \frac \delta 2}  \widetilde {\phi}  (l,y) dy \mbox{ and }C_{**} \lim_{l \rightarrow \infty} \int_{-1 - \frac \delta 2}^{1 + \frac \delta 2}  \widetilde {\phi}  (l,y) dy .$$ By Lemma \ref{lemma_rho}, 
$$ \frac 1 {C_{**} }(2+\delta)\inf_{(x,y) \in \p D_{R0}} \widetilde {\phi} (x,y)  \leq \int_{\p D_{R0} }  x \p_{\nu} \phi ds  \leq C_{**} (2+\delta)\sup_{(x,y) \in \p D_{R0}} \widetilde {\phi} (x,y).$$ Thus, $$  (2+\delta) \frac 1 {C_{*} C_{**}} \sqrt {\epsilon} \leq  \int_{\p D_{R0} }  x \p_{\nu} \phi ds   \leq (2+\delta) C_{*} C_{**} \sqrt {\epsilon}. $$  By the symmetric property,  $\int_{\p D_{L0} }  x \p_{\nu} \phi ds  = \int_{\p D_{R0} }  x \p_{\nu} \phi ds .$ Therefore, the equality \eqref{pro_equal}  implies the desirable result in Proposition \ref{prop_poten_differ_1}. \qed

\section {Proof of Theorem \ref{upp_thm}}  \label{upp_thm_sec}

We begin by defining the domains as  \begin{align*} &{\widetilde \Omega} = \mathbb {R} \times \left(- 1- \frac 1 2 {\delta} , 3+  \frac 3 2 {\delta}\right),~&& {\widetilde \Omega}_4 = \left(-4,4\right) \times \left(- 1- \frac 1 2 {\delta} , 3+  \frac 3 2 {\delta}\right),\\& {\widetilde \Omega}_{4R}= \left(0,4\right) \times \left(- 1- \frac 1 2 {\delta} , 3+  \frac 3 2 {\delta}\right), && {\widetilde \Omega}_{4L}= \left(-4,0\right) \times \left(- 1- \frac 1 2 {\delta} , 3+  \frac 3 2 {\delta}\right)\end{align*} which are used in this proofs. We assume that $H$ is a harmonic function with a periodic gradient as  \eqref{H_condi} and $u$ is a solution to \eqref {main} with the condition \eqref{main_condi}.

\begin{lem} \label{lem3-1} There exists a constant $C$ regardless of $\delta$ and $\epsilon$ such that  \beq 
\left| u \big|_{\p D_{R0}} - u \big|_{\p D_{L0}} \right| = \left| u \big|_{\p D_{R1}} - u \big|_{\p D_{R1}} \right| \leq C \norm {H}_{L^{\infty} (\Omega_4)} \sqrt \epsilon. \notag
\eeq
\end{lem}
\pf 
Let $\widetilde u $ and $r$ be defined by 
$$ \widetilde   u(x,y) = \frac 1 2 \left(u(x,y) + u (x,-y)\right) \mbox{ and } r (x,y) = \frac 1 2 \left(u(x,y) - u (x,-y)\right) $$ and let $$ \widetilde  H (x,y) = \frac 1 2 \left(H(x,y) + H (x,-y)\right). $$ Since $ u = \widetilde  u + r $ and $ r |_{\p D_{L0}} = r |_{\p D_{R0}} = 0 $, Proposition \ref{phi_lemma} implies \begin{align*} u |_{\p D_{R0}} - u |_{\p D_{L0}}  &= \widetilde  u  |_{\p D_{R0}} - \widetilde  u  |_{\p D_{L0}}\\
&= \int_{\partial D_{L0} \cup \partial D_{R0}} {\widetilde  H} \partial_{\nu} \phi  ds\\
&= \int_{\partial D_{L0} \cup \partial D_{R0}} \left({\widetilde  H}(x,y) - {\widetilde  H}(0,0) \right)\partial_{\nu} \phi  ds. \end{align*} It follows from $\p_y {\widetilde  H}(0,0) = 0 $ that  $$ |{\widetilde  H}(x,y)- {\widetilde  H}(0,0) | \leq C \left(  \norm {\nabla H}_{L^{\infty} (\Omega_3)} +  \norm {D^2  H}_{L^{\infty} (\Omega_3)} \right) | x | $$ for any $(x,y) \in \partial D_{L0} \cup \partial D_{R0}$. Thus, 
\begin{align*} \left| u |_{\p D_{R0}} - u |_{\p D_{L0}} \right|  \leq  C \left(  \norm {\nabla H}_{L^{\infty} (\Omega_3)} +  \norm {D^2  H}_{L^{\infty} (\Omega_3)} \right)  \int_{\partial D_{L0} \cup \partial D_{R0}} x  \partial_{\nu} \phi  ds. \end{align*} Here, the standard gradient estimate for harmonic functions implies that 
$$\norm {\nabla H}_{L^{\infty} (\Omega_3)} + \norm {D^2  H}_{L^{\infty} (\Omega_3)}   \leq  C \norm {H}_{L^{\infty} ( (-4,4)\times (-3,3) )}$$  for some $C>0$, since the domain $\Omega_3$ has nonzero distance from $\p ( (-4,4)\times (-3,3) )$, and the periodic property of $\nabla H$ or    \eqref {H=H-H} imply $$ \norm {H}_{L^{\infty} ( (-4,4)\times (-3,3) )} \leq 2 \norm {H}_{L^{\infty} ( \Omega_4)}.$$ Thus, we are done. \qed

\vskip 20pt

The following lemma  provides some maximal principles more general than Lemma \ref {lemma_rho}. 
\begin{lem}  \label{rhoij} Let $\Omega_R$ be as defined in the proof of Proposition \ref{prop_poten_differ_1}. Assume that $\rho_{00}$, $\rho_{10}$, $\rho_{01}$ and  $\rho_{11}$ are harmonic functions defined on $\Omega_R \setminus \overline {D_{R0}} $ with the boundary conditions: $$
 \begin{cases}
 \ds \partial_{\nu}^{i} \rho_{ij} = 0 \quad & \mbox{on } (0,\infty) \times \left\{y~\Big|~ y= -1 - \frac {\delta}  2 ~\mbox{or}~y= 1 + \frac {\delta}  2\right\},  \\ \nm \ds
\partial_{\nu}^{j}  \rho_{ij} = 0   \quad &\mbox{on } \left\{0\right\}  \times \left(-1 - \frac {\delta}  2,~ 1 + \frac {\delta}  2 \right) ,\end{cases}
 $$ and also satisfying 
$$ \int_{\Omega_R \setminus \overline {D_{R0}}} \left|\nabla \rho_{ij} \right |^2  dxdy < \infty $$ for any $i,~j =0,~1$, where $\p_{\nu} ^0 u = u$ and $\p_{\nu} ^1 u = \p_{\nu} u$. Then, $$ \norm{ \rho_{ij} }_{L^{\infty}(\Omega_R \setminus \overline {D_{R0}})}   \leq  \norm{ \rho_{ij} }_{L^{\infty}(\p D_{R0})} .$$ 
\end{lem}  This lemma can also  be derived easily by  the same function $\varPhi_R: B_1^+ (0,0) \rightarrow \Omega_R$ used in Lemmas \ref{lemma_phi_R} and  \ref {lemma_rho}, the maximum principle and the Hopf Lemma. The proof of this lemma is left as an exercise for the reader.

\begin{lem} \label{lem3-2}   Let $a_*$ be the constant defined as  $$a_* = u \left(-\frac \epsilon 2, 0\right)-H \left(-\frac \epsilon 2, 0\right)= u \big|_{\p D_{L0}} - H  \left(-\frac \epsilon 2, 0\right).$$ Then,  \begin{align*}
\norm{ u - H- a_*}_{L^{\infty} \left(\widetilde {\Omega} \setminus \overline{ D_{L0} \cup D_{L1} \cup D_{R0} \cup D_{R1}  }\right) } \leq C \norm {H}_{L^{\infty} ({\Omega}_4)}
\end{align*} and 
$$\norm{ u - u \big|_{\p D_{L0}} }_{L^{\infty} \left(\widetilde {\Omega}_4 \setminus \overline{ D_{L0} \cup D_{L1} \cup D_{R0} \cup D_{R1}  }\right) } \leq C \norm {H}_{L^{\infty} ({\Omega}_4)}.$$
\end{lem}
\pf We define the notations $(\cdot)_e$ and $(\cdot)_o$ as follows:
$$ (v)_e (x,y) = \frac 1 2  \left(v(x,y) + v(x,-y)\right) $$ and
$$ (v)_o (x,y) = \frac 1 2  \left(v(x,y) - v(x,-y)\right) $$ for a function $v$ defined in $ \overline {\Omega} \setminus \left({ D_{L0} \cup D_{R0}}\right)$.
Then, $$ u - H = (u-H)_e +  (u-H)_o.$$ 

\par First, we estimate $\norm {(u-H)_e - a_*}_{L^{\infty} (  \Omega \setminus \overline {D_{L0} \cup  D_{R0}})} $.  For any $(x,y) \in \p D_{L0}$,
\begin{align*}
&(u-H)_e (x,y) - a_* \\& = (u-H)_e (x,y) - u \big|_{\p D_{L0}} + H  \left(-\frac \epsilon 2, 0\right)\\& = (u)_e (x,y) - u \big|_{\p D_{L0}} -(H)_e (x,y) + (H)_e  \left(-\frac \epsilon 2, 0\right)
\\& = -(H)_e (x,y) + (H)_e  \left(-\frac \epsilon 2, 0\right),
\end{align*} and for any $(x,y) \in \p D_{R0}$,
\begin{align*}
&(u-H)_e (x,y) - a_* \\& = (u-H)_e (x,y) - u \big|_{\p D_{L0}} + H  \left(-\frac \epsilon 2, 0\right)\\& = u \big|_{\p D_{R0}} - u \big|_{\p D_{L0}} -(H)_e (x,y) + (H)_e  \left(-\frac \epsilon 2, 0\right).
\end{align*} Thus, the equality \eqref{pro_equal} implies that 
\beq \norm {(u-H)_e - a_*}_{L^{\infty} (\p D_{L0} \cup \p D_{R0})} \leq 4\norm {H}_{L^{\infty} (\p D_{L0} \cup \p D_{R0})} . \label{djWjekakwncls}\eeq
Note that   $(u-H)_e - a_*$ is a harmonic function in $\Omega \setminus \overline { D_{L0} \cup  D_{R0}}$ with 
$$
 \begin{cases}
 \ds \partial_{\nu} \left((u-H)_e - a_* \right) = 0 & \mbox{on } (-\infty,\infty) \times \left\{y~\Big|~ y= -1 - \frac {\delta}  2 ~\mbox{or}~y= 1 + \frac {\delta}  2\right\},  \\ \nm \ds\int_{\Omega_R \setminus \overline {D_{L0}\cup D_{R0}}} \left|\nabla \left((u-H)_e - a_* \right) \right |^2  &dxdy < \infty \end{cases}
 $$ due to a periodic property of $\nabla u$ and $\nabla H$. Applying $(u-H)_e - a_*$ to $\rho_{10}+\rho_{11}$ in Lemma \ref{rhoij}, the bound \eqref{djWjekakwncls} implies $$ \norm {(u-H)_e - a_*}_{L^{\infty} (  \Omega \setminus \overline {D_{L0} \cup  D_{R0}})} \leq 4\norm {H}_{L^{\infty} (\p D_{L0} \cup \p D_{R0})} .$$
 
 \vskip 10pt
 
 \par Second, we estimate $\norm {(u-H)_o }_{L^{\infty} (  \Omega \setminus \overline {D_{L0} \cup  D_{R0}})} $.  On $ \p D_{L0} \cup \p D_{R0} $, $$ (u-H)_o  = -(H)_o  ,$$ since $u$ is constant on $\p D_{L0} $  and $\p D_{R0}$, respectively.  Meanwhile, the equality \eqref{proof_first_prop} in the proof of Proposition \ref{prop_poten_differ_2} means that the harmonic function $(u-H)_o $ satisfies 
$$ (u-H)_o (x,y) =0 \mbox{ on } (-\infty,\infty) \times \left\{y~\Big|~ y= -1 - \frac {\delta}  2 ~\mbox{or}~y= 1 + \frac {\delta}  2\right\}.$$ Since $  (u-H)_o = -(H)_o $ on $\p D_{L0} \cup \p D_{R0}$ and $\int_{\Omega_R \setminus \overline {D_{L0}\cup D_{R0}}} \left|\nabla \left((u-H)_e - a_* \right) \right |^2 dxdy < \infty $, the results on $\rho_{00}$ and  $\rho_{01}$ in Lemma \ref{rhoij} yield 
$$ \norm {(u-H)_o}_{L^{\infty} (  \Omega \setminus \overline {D_{L0} \cup  D_{R0}})} \leq \norm {(u-H)_o}_{L^{\infty} (   {\p D_{L0} \cup  \p D_{R0}})} \leq  \norm {H}_{L^{\infty} (\p D_{L0} \cup \p D_{R0})} .$$

\par  Combining the first and second cases,  $$ \norm {u-H-a_*}_{L^{\infty} (  \Omega \setminus \overline {D_{L0} \cup  D_{R0}})} \leq  5 \norm {H}_{L^{\infty} (\p D_{L0} \cup \p D_{R0})} .$$ The first inequality in this lemma can be derived by  \eqref{proof_first_prop} as follows: \begin{align*} &\norm{ u - H- a_*}_{L^{\infty} \left(\widetilde {\Omega} \setminus \overline{ D_{L0} \cup D_{L1} \cup D_{R0} \cup D_{R1}  }\right) }\\&=  \norm {u-H-a_*}_{L^{\infty} (  \Omega \setminus \overline {D_{L0} \cup  D_{R0}})} \\ &\leq  5 \norm {H}_{L^{\infty} (\p D_{L0} \cup \p D_{R0})}. \end{align*} The second inequality also follows immediately  so that \begin{align*}  &\norm{ u - u|_{\p D_{L0}}}_{L^{\infty} \left(\widetilde {\Omega}_4 \setminus \overline{ D_{L0} \cup D_{L1} \cup D_{R0} \cup D_{R1}  }\right) }\\ &\leq \norm{ u - H- a_*}_{L^{\infty} \left(\widetilde {\Omega} \setminus \overline{ D_{L0} \cup D_{L1} \cup D_{R0} \cup D_{R1}  }\right) }+ \norm{ H - H \left( -\frac \epsilon 2, 0   \right)}_{L^{\infty} \left(\widetilde {\Omega}_4 \setminus \overline{ D_{L0} \cup D_{L1} \cup D_{R0} \cup D_{R1}  }\right) }  \\ &\leq  5 \norm {H}_{L^{\infty} (\p D_{L0} \cup \p D_{R0})} +  2 \norm{ H }_{L^{\infty} \left(\widetilde {\Omega}_4 \right) }\\ &\leq  7 \norm{ H }_{L^{\infty} \left(\widetilde {\Omega}_4 \right) }\\ &\leq  14 \norm{ H }_{L^{\infty} \left( {\Omega}_4 \right) } ,\end{align*} due to \eqref {H=H-H}. We are done. \qed

\subsection {The proof of \eqref{upp_thm_1st_inequ}}
The potential difference $u |_{D_{R 1}} - u |_{D_{R0}} $ was evaluated exactly in Proposition \ref{prop_poten_differ_2}. The value has very different nature from the cases of finite number of inclusions, and also results in much stronger concentration than finite cases. In this proof, we establish an asymptote of $\nabla u$  from the potential difference.  Indeed, a nice method to get an asymptote was already introduced by Kang, Lim, Yun in the case of two circular inclusions  in \cite{KLY}, and Bao, Li, Yin in \cite{BLY} showed the boundedness of the gradient in the case of no potential difference. In this proof, we modify these methods to apply to our problem, and obtain an asymptote describing the stronger concentration. Hence, the potential difference evaluated in Proposition \ref{prop_poten_differ_2} plays the most important role in the result.

\par To establish the asymptote, we consider the decomposition of $\nabla u$ into two terms as 
$$ \nabla u = \alpha_h \nabla \phi_h +  \nabla u_h, $$ where $\alpha_h $, $ \phi_h$ and $ u_h$ are definded below.  The function $ \phi_h $ has a high concentration in between $D_{R0}$ and $D_{R1}$, and is also easy to handle. In this proof, we estimate the coefficient $\alpha_h$ and show that  $\nabla u_h$ is bounded regardless of $\epsilon$ and $ \delta$. Thus, we can establish the desirable asymptote  \eqref{upp_thm_1st_inequ}.

\par We define  $\alpha_h $, $ \phi_h$ and $ u_h$ and set the decomposition up. Let $\phi_h (x,y)$ be the unique solution to  \beq
 \begin{cases}
 \ds \Delta \phi_h  = 0 \quad  & \mbox{in } \mathbb{R}^2 \setminus \overline{D_{R0} \cup D_{R1}}\\
 \nm \ds  \phi_h = \mbox{ a constant}   \quad & \mbox{on } \partial  D_{R1}, \\ \nm  \ds \phi_h = - \phi_h \Big|_{ \partial  D_{R1}}   ~ \quad & \mbox{on } \partial  D_{R0},   \\ \nm \ds \int_{\partial D_{R1}} \partial_{\nu} \phi_h  ds =  -\int_{\partial D_{R0}} \partial_{\nu}  \phi_h  ds &=  \frac {2\pi} {\sqrt \delta}, \\  \nm \ds   \phi_h ({\bf x}) = O\left({\frac 1 {|{\bf x}|}}\right) & \mbox{as } {|\bf x|} \rightarrow \infty. \end{cases}\eeq in the same as \eqref {eqn_phi_n_before} and \eqref{eqn_phi_n}. The solution can be expressed as  $$\phi_h (x,y)=  \frac 1 {\sqrt \delta} \left( \log \left|(x,y) - \left(1+\frac \epsilon 2 , 1+\frac \delta 2 - p_h \right)\right| -\log \left|(x,y) - \left(1+\frac \epsilon 2 , 1+\frac \delta 2 + p_h \right)\right| \right)$$ where 
$$ p_h = \sqrt \delta + O (\delta)$$ for small $\delta$. Let $\alpha_h$ be  the constant as 
$$\alpha_h = \frac {u|_{\p D_{R1}} - u|_{\p D_{R0}} }{\phi_h|_{\p D_{R1}} - \phi_h|_{\p D_{R0}}  } ,$$  and we define a harmonic function $u_{h}$ as $u_{h} = u - \alpha_h \phi_h  - (u - \alpha_h \phi_h)|_{\p D_{R0}}.$ The solution $u$ is decomposed into $ \alpha_h \phi_h + (u- \alpha_h \phi_h)|_{\p D_{R0}}$ and $u_{h}$ as follows:
$$ u = \left(\alpha_h \phi_h + (u- \alpha_h \phi_h)|_{\p D_{R0}}\right)  + u_h.$$ Hence, $$\nabla  u = \alpha_h \nabla  \phi_h    + \nabla u_h.$$

\par From the defintion of $\alpha_h$, two functions $ \alpha_h \phi_h$ and $u$ have the same potential difference between $\p D_{R1}$ and $\p D_{R0}$, and $u_h$ has no difference between the boundaries so that $ \alpha_h \phi_h |_{\p D_{R1}}    - \alpha_h \phi_h |_{\p D_{R0}}=u|_{\p D_{R1}}    - u|_{\p D_{R0}}$ and $ u_h |_{\p D_{R1}} -  u_h |_{\p D_{R0}} = 0$. Indeed, this means that $\nabla u$ is dominated by $\alpha_h \nabla \phi_h$.  By direct calculation and the definition of $N_h$, there is a constant $C$ regardless of $\epsilon$ and $\delta$ such that $$  \left|  \nabla \phi_h - 2\left(-2 ~  \frac { \left(x-1-\frac \epsilon 2 \right)\left(y-1-\frac \delta 2 \right)}{\left(\left(x-1-\frac \epsilon 2 \right)^2 + \delta\right)^2} , \frac 1 { \delta +\left(x -1 - \frac {\epsilon} 2 \right)^2} \right) \right| \leq C $$ and  \begin {align*}&\left| 2\sqrt {\delta}~  \frac { \left(x-1-\frac \epsilon 2 \right)\left(y-1-\frac \delta 2 \right)}{\left(\left(x-1-\frac \epsilon 2 \right)^2 + \delta\right)^2}    \right|=|S(x,y)| \\& \leq   \left|   \frac { y-1-\frac \delta 2}{\left(x-1-\frac \epsilon 2 \right)^2 + \delta}    \right|  \leq  2 \left|   \frac { \left(x-1-\frac \epsilon 2 \right)^2 }{\left(x-1-\frac \epsilon 2 \right)^2 + \delta}    \right|  \leq 2 \end{align*} in $N_h$.  By Proposition \ref{prop_poten_differ_2},  direct calculation of $\phi_h$ implies $$\left|\alpha_h - \frac 1 2 \left(H (0,1) - H (0,-1)\right) \right| \leq C  \delta \norm {\nabla H}_{L^{\infty} ({\Omega}_3)}.$$ By \eqref{H=H-H}, a  standard gradient estimate for harmonic functions yields  $\norm {\nabla H}_{L^{\infty} ({\Omega}_3)}  \leq  C \norm {H}_{L^{\infty} (\widetilde{\Omega}_4)}  \leq 3 C\norm { H}_{L^{\infty} ({\Omega}_4)} $ in the same way as the proof of Lemma \ref{lem3-1}. Hence, \beq \left|\alpha_h - \frac 1 2 \left(H (0,1) - H (0,-1)\right) \right| \leq C \delta \norm { H}_{L^{\infty} ({\Omega}_4)}. \label{gunpla_choishoonsil}\eeq

\par The remainder of the proof is deducated only to prove the boundedness of $\nabla u_h$ such that $$\norm{\nabla u_h}_{L^\infty (N_h )} \leq  C \norm {H}_{L^{\infty} ({\Omega}_4)}$$ for some constant $C$. Then, we obtain the desirable  \eqref{upp_thm_1st_inequ}.   Some propeties of $u_h$ are considered before proving the boundedness. From the defintion of $\alpha_h$, $$ u_h |_{\p D_{R1}} -  u_h |_{\p D_{R0}} = 0,$$ $$ \alpha_h \phi_h |_{\p D_{R1}}    - \alpha_h \phi_h |_{\p D_{R0}}=u|_{\p D_{R1}}    - u|_{\p D_{R0}}$$ and  $|\alpha_h | \leq 2 \norm {H}_{L^{\infty} ({\Omega}_4)}$ by \eqref{gunpla_choishoonsil}.  Lemma \ref{lem3-2} implies that \begin{align} &\norm { u_{h }  }_{L^{\infty} \left( \widetilde{ \Omega}_{4R }\setminus  \overline {D_{R0} \cup  D_{R1} } \right) } \notag\\& \leq  \norm { u  - u|_{\p D_{R0}} }_{L^{\infty} \left( \widetilde{ \Omega}_{4R }\setminus  \overline {D_{R0} \cup  D_{R1} } \right)}+ \left| \alpha_h \phi_h \big|_{\p D_{R0}}\right| + \norm { \alpha_h \phi_h  }_{L^{\infty} \left( \widetilde{ \Omega}_{4R} \setminus  \overline {D_{R0} \cup  D_{R1} } \right)}  \notag \\& \leq C \norm {H}_{L^{\infty} ({\Omega}_{4})}.\label {3-1independent} \end{align}  From defintion, \beq  u_{h} = 0 \mbox{ on } \p D_{R0} \cup \p D_{R1}. \label {3-2independent} \eeq

\par Dealing with the boundedness of $\nabla u_h$, we decompose $ u_{h }$  into two functions $u_{+} $ and $u_{- }$ as 
$$  u_{h }  =  u_{+} +  u_{- } , $$ where $u_+$ and $u_-$ are the harmonic functions given as  \beq
 \begin{cases} \ds \triangle  u_{+} =\triangle u_{- } =0& \mbox { in }   \widetilde {\Omega}_{4R} \setminus \overline {D_{R0} \cup D_{R1}}, \\  \ds u_{+}\big |_{\p D_{R0} \cup \p D_{R1}}= u_{-}\big |_{\p D_{R0} \cup \p D_{R1}} = 0,& \\ \ds u_{+}(x,y) =  \max \{u_{h}(x,y),0 \}\geq 0~&\mbox{for any }(x,y)\in \p \widetilde {\Omega}_{4R} , \\ \ds  u_{-}(x,y) =  \min \{u_{h}(x,y),0 \}\leq 0 ~&\mbox{for any }(x,y)\in \p \widetilde {\Omega}_{4R}.\end{cases} \notag
 \eeq  Then, $$ u_{+} \geq 0  \mbox{ and }  u_{-} \leq 0 \mbox { in } \widetilde {\Omega}_{4R} \setminus \overline {D_{R0} \cup D_{R1}}.$$

\par In order to derive the boundedness of $\nabla u_+$ in $N_h$, we estimate $\nabla u_+$  on the boundary $\p N_h$ which consists of four curves $\p D_{R1} \cap \p N_h$, $\p D_{R0} \cap \p N_h$, $\left\{1+ \frac {\sqrt 3 }{2} + \frac \epsilon 2 \right\}\times [\frac 1 2 , \frac 3 2 +\delta]$ and $\left\{ 1- \frac {\sqrt 3 }{2} + \frac \epsilon 2\right\}\times [\frac 1 2 , \frac 3 2 +\delta]$. We use $ u_{+0} $ and $u_{+1}$ defined as
\beq
 \begin{cases}
 \ds  \triangle u_{+0} = 0 &\mbox{ in } \widetilde {\Omega}_{4R} \setminus \overline {D_{R0}},  \\ \nm \ds  u_{+0} = u_+  &\mbox{on }\p \widetilde {\Omega}_{4R} ,  \\ \nm \ds  u_{+0} = 0   &\mbox{on } \p D_{R0},  \end{cases} \notag
\mbox{ and }
 \begin{cases}
 \ds  \triangle u_{+1} = 0 &\mbox{ in } \widetilde {\Omega}_{4R} \setminus \overline {D_{R1}},  \\ \nm \ds  u_{+1} = u_+  &\mbox{on }\p \widetilde {\Omega}_{4R} ,  \\ \nm \ds  u_{+1} = 0   &\mbox{on } \p D_{R1}. \end{cases} \notag
 \eeq It follows from definitions and \eqref {3-1independent} that
\begin{align}  &\norm{ u_{+0}}_{L^{\infty }(\widetilde {\Omega}_{4R} \setminus \overline {D_{R0}})}  +  \norm{ u_{+1}}_{L^{\infty }(\widetilde {\Omega}_{4R} \setminus \overline {D_{R1}})} \notag \\& \leq 2  \norm{ u_{h}}_{L^{\infty }(\widetilde {\Omega}_{4R} \setminus \overline {D_{R0} \cup D_{R1}})} \leq  C \norm{H}_{L^{\infty }( {\Omega}_{4})} . \label{3-4-asdfg} \end{align} Since $u_{+0} -  u_+ = 0 $ on  $\p \widetilde {\Omega}_{4R} \cup \p D_{R0}$  and $u_{+0} -  u_+ \geq  0 $ on $ \p D_{R1}$,  \beq 0\leq  u_+  \leq u_{+0} ~~~ \mbox{ in }\widetilde {\Omega}_{4R} \setminus \overline {D_{R0} \cup D_{R1}}.\label{911math}\eeq  Since $ u_{+0} - u_+  = u_+ =0$ on $\p D_{R0}$, the functions $u_{+0} -  u_+$ and $ u_+$ attain  the minimal value $0$ on $\p D_{R0}$. The Hopf's lemma thus implies that 
 $ 0 \geq \p_{\nu} u_{+} \geq   \p_{\nu} u_{+0}$. Thus, \beq  0 \leq |\nabla u_{+} |\leq   |\p_{\nu} u_{+0}|~\mbox { on } \p D_{R0}, \label{u_++_bound}\eeq and similarly \beq 0 \leq |\nabla u_{+} |\leq   |\p_{\nu} u_{+1}| ~\mbox { on } \p D_{R1}.\label{u_--_bound}\eeq 
 Since $u_{+0}=0$ on $\p  D_{R0}$ and $u_{+1}=0$ on $\p  D_{R1}$, the Kelvin transform can extend the functions $u_{+0}$, and $u_{+1}$, into harmonic functions $\widetilde {u}_{+0}$, and $\widetilde {u}_{+1}$, defined open sets containing $\p D_{R0}$, and  $\p D_{R1}$, respectively. For any $ (x_0, y_0) \in \p D_{R0} \cap \p N_h$, the extended function $\tilde{u}_{+0}$ is defined in $B_{\frac 1 8} (x_0,y_0)$. A gradient estimate for harmonic functions and \eqref {3-4-asdfg} yield \begin{align*} &\left |\nabla  {u}_{+0} (x_0,y_0)\right |  =\left |\nabla \widetilde {u}_{+0} (x_0,y_0)\right | \leq C_1 \left(\frac 1 8 \right)^{-1} \sup_{B_{\frac 1 8} (x_0,y_0)} \left| \tilde u_{+0} (x,y) \right|  \\ & \leq C_2 \norm{ u_{+0}}_{L^{\infty }(\widetilde {\Omega}_{4R} \setminus \overline {D_{R0}})}   \leq C_3 \norm{H}_{L^{\infty }( {\Omega}_{4})}. \end{align*} Thus, \eqref {u_++_bound} implies \beq \norm{ \nabla u_{+} }_{L^{\infty} (\p D_{R0} \cap \p N_h)} \leq  C  \norm{H}_{L^{\infty }( {\Omega}_{4})}, \label{enrowndgksk} \eeq and in the same way, \eqref {3-4-asdfg} and \eqref{u_--_bound} yield \beq \norm{ \nabla u_{+} }_{L^{\infty} (\p D_{R1} \cap \p N_h)} \leq  C  \norm{H}_{L^{\infty }( {\Omega}_{4})}.\label{enrowndenqjsWo} \eeq  Hence, we have the upper bounds for $ |\nabla u_{+} |$ on each boundaries  $\p D_{R1} \cap \p N_h$ and $\p D_{R0} \cap \p N_h$ as above. Meanwhile, we estimate $ |\nabla u_{+} |$ on two vertical line segments $\left\{1+ \frac {\sqrt 3 }{2} + \frac \epsilon 2 \right\}\times [\frac 1 2 , \frac 3 2 +\delta]$ and $\left\{ 1- \frac {\sqrt 3 }{2} + \frac \epsilon 2\right\}\times [\frac 1 2 , \frac 3 2 +\delta]$ which are the remainder boundaries of $\p N_h$. Since $u_+ = 0$ on $\p D_{R0} \cup \p D_{R1}$, the Kelvin transform extends $u_+$ to a harmonic function $\tilde {u}_+$  defined in an open set containing $\p {\left(D_{R0} \cup D_{R1} \right)}$ as well as  $ \widetilde {\Omega}_{4R}  \setminus  \overline{\left(D_{R0} \cup D_{R1} \right)}$. For any point $(x_0,y_0) $ on the vertical line segments above, the extended harmonic function $\tilde u_+$ is defined in the open disk $B_{\frac 1 8 } (x_0,y_0)$. A gradient estimate for harmonic functions and \eqref{911math} thus yield \begin{align}  &|\nabla u_+ (x_0,y_0) | =|\nabla \tilde u_+ (x_0,y_0) | \leq C_1 \left(\frac 1 {8}\right)^{-1} \sup_{B_{\frac 1 8} (x_0,y_0)} \left| \tilde u_+ (x,y) \right|   \notag  \\&\leq C_2 \norm{u_{+}}_{L^\infty (\widetilde {\Omega}_{4R} \setminus \overline {D_{R0}\cup D_{R1}})}  \leq C_2   \norm{u_{+0}}_{L^\infty (\widetilde {\Omega}_{4R} \setminus \overline {D_{R0}})}. \label{3-5_+_+_+independent} \end{align} By the definitions of $u_{+0}$ and  $u_{+1}$, and by \eqref{3-1independent}, 
$\norm{u_{+0}}_{L^\infty (\widetilde {\Omega}_{4R} \setminus \overline {D_{R0}})} + \norm{u_{+1}}_{L^\infty (\widetilde {\Omega}_{4R} \setminus \overline {D_{R1}})} \leq C_3  \norm {H}_{L^{\infty} ({\Omega}_{4})} $. Hence, \eqref{enrowndgksk}, \eqref{enrowndenqjsWo} and \eqref{3-5_+_+_+independent} result in a gradient estimate on the boundary $\p N_h$ as  $$ \norm{\nabla u_{+}}_{L^\infty (\p N_h )}  \leq C_4  \norm {H}_{L^{\infty} ({\Omega}_{4})}.$$ By the maximal principle, 
$$ \norm{\nabla u_{+}}_{L^\infty ( N_h )}  \leq C_4  \norm {H}_{L^{\infty} ({\Omega}_{4})}.$$ In the same way,  we also get  $$\norm{\nabla u_{-}}_{L^\infty (N_h )} \leq  C_5 \norm {H}_{L^{\infty} ({\Omega}_4)}.$$ Therefore, we obtain $$\norm{\nabla u_h}_{L^\infty (N_h )} \leq  C_6 \norm {H}_{L^{\infty} ({\Omega}_4)}.$$ We are done. \qed

\subsection {The proof of \eqref{upp_thm_2nd_inequ}}
An estimate for the potential difference $u|_{\p D_{R0}} - u|_{\p D_{L0}}$ was obtained in Lemma \ref{lem3-1}. We repeat the same method as the proof of \eqref{upp_thm_1st_inequ} to establish the asymptote from the potential difference. Thus, this proof also begins at the decomposition as $$ u = \beta_v \nabla  \phi_v    + \nabla u_v.$$ Here, $\phi_v$ is  the unique solution to \beq
 \begin{cases}
 \ds \Delta \phi_v  = 0 \quad  & \mbox{in } \mathbb{R}^2 \setminus \overline{D_{L0} \cup D_{R0}}\\
 \nm \ds  \phi_v = \mbox{ a constant}   \quad & \mbox{on } \partial  D_{R0}, \\ \nm  \ds \phi_v = - \phi_v \Big|_{ \partial  D_{R0}}   ~ \quad & \mbox{on } \partial  D_{L0},   \\ \nm \ds \int_{\partial D_{R0}} \partial_{\nu} \phi_v  ds =  -\int_{\partial D_{L0}} \partial_{\nu}  \phi_h  ds &=  {2\pi}, \\  \nm \ds   \phi_v ({\bf x}) = O\left({\frac 1 {|{\bf x}|}}\right) & \mbox{as } {|\bf x|} \rightarrow \infty. \end{cases}\eeq Then, $$\phi_v (x,y)=  \log \left|(x,y) + \left(p_v,0 \right)\right| -\log \left|(x,y) - \left(p_v,0 \right)\right| $$ where 
$$ p_v = \sqrt \epsilon + O (\epsilon).$$  Let $\beta_v$ be  the constant as 
$$\beta_v = \frac {u|_{\p D_{R0}} - u|_{\p D_{L0}} }{\phi_v|_{\p D_{R0}} - \phi_v|_{\p D_{L0}}  } $$  and we define a harmonic function $u_{v}$ as $$u_{v} = u - \beta_v \phi_v  - (u - \beta \phi_v)|_{\p D_{R0}}.$$ The solution $u$ is decomposed into $ \beta_v \phi_v + (u- \beta_v \phi_v)|_{\p D_{R0}}$ and $u_{v}$ as 
$$ u = \left(\beta_v \phi_v + (u- \beta_v \phi_v)|_{\p D_{R0}}\right)  + u_v.$$ Hence, we obtain the desirable decomposition  $$ u = \beta_v \nabla  \phi_v    + \nabla u_v.$$

\par  By direct calculation, there is a constant $C$ regardless of $\epsilon$ and $\delta_v$ such that $$  \left|   \nabla  \phi_v - 2 \frac {\sqrt{\epsilon}} { \epsilon +y^2}   \left(1,0\right) \right| \leq C ,$$ and Lemma \ref{lem3-1} implies $|\beta | \leq 3 \norm {H}_{L^{\infty} ({\Omega}_4)}$. If we prove the boundedness of $\nabla u_v$ such that  $$\norm{\nabla u_v}_{L^\infty (N_v )} \leq  C \norm {H}_{L^{\infty} ({\Omega}_4)}$$ for some constant $C$, then we can obtain the main result \eqref{upp_thm_2nd_inequ}.

\par The remainder of the proof is  to prove the boundedness of $\nabla u_v$.  From the defintion of $\beta_v$, we have  $ \beta_v \phi_v |_{\p D_{R0}}    - \beta_v \phi_v |_{\p D_{L0}}=u|_{\p D_{R0}}    - u|_{\p D_{L0}}$, $ u_v |_{\p D_{R0}} -  u_v |_{\p D_{L0}} = 0$ and  $|\beta_v | \leq 3 \norm {H}_{L^{\infty} ({\Omega}_4)}$ by  Lemma \ref{lem3-1}.  Similarly to \eqref  {3-1independent},  Lemma \ref{lem3-2} implies that \begin{align} &\norm { u_{v}  }_{L^{\infty} \left({ \Omega}_{4}\setminus  \overline {D_{R0} \cup  D_{L0} } \right) }  \leq C \norm {H}_{L^{\infty} ({\Omega}_{4})}.\label {3-9independent} \end{align}  From defintion, \beq  u_{v} = 0 \mbox{ on } \p D_{L0} \cup \p D_{R0}. \label {3-10independent} \eeq They are the conditions analogous to \eqref{3-1independent} and \eqref{3-2independent}. In the same way as the proof of \eqref{upp_thm_1st_inequ}, we have $$\norm{\nabla u_v}_{L^\infty (N_v )} \leq  C \norm {H}_{L^{\infty} ({\Omega}_4)}.$$ We are done. \qed

\section {Proofs of Theorems \ref{low_thm}} \label {low_thm_sec}
The proof is mainly concerned with the second equality \eqref {low_thm_2nd_eq} in $N_v$, since the first equality \eqref {low_thm_1st_eq} in $N_h$ follows immediately from Theorem \ref{upp_thm}. Owing to the linearity of problem, we consider two cases when $H(x,y)= x$, and when $H(x,y)= y$, separately. Let $u_{a}$ and $u_b$ be the solutions for $H(x,y)=x$ and $H(x,y)=y$, respectively.

\par In the first case when $H(x,y)= x$ for $(x,y)\in \mathbb{R}^2$,  Theorem \ref {upp_thm} presents a constant $\mu_0$ satisfying 
\beq \nabla u_a (x,y) =  \mu_0 \frac {\sqrt \epsilon}{ \epsilon + y^2 }   \left(1,0\right) + R_{a2} (x,y) \label{upp_thm_wndjdhsrjt} \eeq for $(x,y)\in N_v$, while $\norm {R_{a2} (x,y)}_{L^{\infty}(N_v)} $ is bounded regardless of small $\epsilon>0$ and $\delta>0$. Proposition  \ref{prop_poten_differ_1} provides a positive constant $C_1$ regardless of $\epsilon$ and $\delta$ such that $$ \frac 1 {C_1} \sqrt \epsilon \leq   u_a |_{D_{R 0}} - u_a |_{D_{L 0}} \leq C_1 \sqrt \epsilon.$$ By the mean value theorem, there exists a point $(x_a,0)\in N_v$ such that $-\frac 1 2 \epsilon< x_a< \frac 1 2 \epsilon$ and $$ \frac 1 {C_1} \frac 1 {\sqrt \epsilon} \leq \p_x u_a(x_a,0) \leq C_1 \frac 1 {\sqrt \epsilon}.$$ By  \eqref{upp_thm_wndjdhsrjt}, the coefficient $\mu_0$ is bounded below as   $$  \mu_0\geq   \frac 1 {C_1} -  \sqrt {\epsilon} ~\norm {R_{a2} (x,y)}_{L^{\infty}(N_v)}    \geq  \frac 1 2 \frac 1 {C_1} $$ for small $\epsilon>0$ due to the boundedness of  $\norm {R_{a2} (x,y)}_{L^{\infty}(N_v)} $. Theorem \ref {upp_thm} provides an upper bound for $\mu_a$ so to obtain a constant $C_2>0$  satisfying $$ \frac 1 {C_2} \leq \mu_0 \leq C_2 $$ regardless of $\epsilon$ and $\delta$. Hence, we have the estimate \eqref {low_thm_2nd_eq} for $\nabla u_a$ in $N_v$ with \eqref{low_thm_3rd_eq}.

\par In the second case when $H(x,y)=  y$ for $(x,y)\in \mathbb{R}^2$, it follows from Theorem \ref {upp_thm}  that \beq \nabla  u_b (x,y) =   {\mu}_{b}  \frac {\sqrt \epsilon}{ \epsilon + y^2 } (1,0)   +  {R_{b2} }(x,y) ~\mbox{for any } (x,y)\in N_v \label{last_thm_asympto_aaa}\eeq  and for a proper constant constant $ {\mu}_{b}$, and $\norm{  {R_{b2}}}_{L^{\infty}(N_v) }$ is bounded regardless of $\epsilon$ and $\delta$. By Proposition \ref {phi_lemma} for $H=y$, 
$$ \frac 1 {\sqrt \epsilon}\int_{-\frac \epsilon 2}^{\frac \epsilon 2}  \p_x u_{b}  (x,0) dx  =  \frac 1 {\sqrt \epsilon} \left( u_{b} |_{\p D_{R0 }} -u_{b} |_{\p D_{L0}} \right)= 0.$$ Applying \eqref{last_thm_asympto_aaa} to here, 
$$ \left | {\mu}_{b} \right| \leq  C_3 \sqrt \epsilon.$$ Applying the inequality to \eqref{last_thm_asympto_aaa},  there exists a constant $C_4$  regardless of $\epsilon$ and $\delta$ such  that $$\norm{\nabla u_{b}}_{L^{\infty}(N_{v})} \leq C_4.$$ 

\par Therefore, in the case when $H(x,y) = a x + b y$ in $\mathbb{R}^2$, we have the desirable asymptote as  \begin{align*}\nabla u &= a \nabla u_a + b \nabla u_b \\&= a  \mu_0 \frac {\sqrt \epsilon}{ \epsilon + y^2 }   \left(1,0\right)  + R_2  \end{align*}  in $N_v$ and the remainder term $R_2$ is bounded, since  $$ R_2 =  a  R_{a2} + b \nabla u_b + b R_{b2}.$$
\qed

\section*{Acknowledgement}

This research was supported
by Basic Science Research Program through the National Research Foundation of Korea (NRF)
funded by the Ministry of Education \\ NRF-2015R1D1A1A01059212


\end{document}